\documentclass[12pt,onecolumn]{IEEEtran}
\usepackage{amsmath,amssymb,yfonts,psfrag,latexsym,dsfont,graphicx,bbm} 

\begin{document}

\newtheorem{thm}{Theorem}
\newtheorem{cor}[thm]{Corollary}
\newtheorem{lemma}[thm]{Lemma}
\newtheorem{prop}[thm]{Proposition}
\newtheorem{problem}[thm]{Problem}
\newtheorem{remark}[thm]{Remark}
\newtheorem{defn}[thm]{Definition}
\newtheorem{ex}[thm]{Example}

\newcommand{\Rnp}{{\mathbb R}^{n+1}}
\newcommand{\mR}{{\mathbb R}}
\newcommand{\mZ}{{\mathbb Z}}
\newcommand{\mS}{{\mathbb S}}
\newcommand{\mI}{{\mathbb I}}
\newcommand{\mC}{{\mathbb C}}
\newcommand{\fR}{{\mathfrak R}}
\newcommand{\fG}{{\mathfrak G}}
\newcommand{\fT}{{\mathfrak T}}
\newcommand{\fC}{{\mathfrak C}}
\newcommand{\fF}{{\mathfrak F}}
\newcommand{\fM}{{\mathfrak M}}
\newcommand{\cM}{{\mathfrak M}}
\newcommand{\cS}{{\cal S}}
\newcommand{\cR}{{\cal R}}
\newcommand{\cU}{{\cal U}}
\newcommand{\cG}{{\cal G}}
\newcommand{\cC}{{\cal C}}
 \newcommand{\cN}{{\cal N}}
 \newcommand{\cK}{{\cal K}}
 \newcommand{\cI}{{\cal I}}
 \newcommand{\rleft}{{\rm left}}
 \newcommand{\rright}{{\rm right}}
 \newcommand{\trace}{{\rm trace}}
 \newcommand{\symmetric}{{\rm Herm}}
\newcommand{\Real}{{\Re}e\,}

\title{Relative Entropy and the multi-variable multi-dimensional Moment Problem}
\author{Tryphon T. Georgiou \thanks{Department of Electrical and Computer Engineering, University of Minnesota, Minneapolis, MN 55455; {\tt tryphon@ece.umn.edu}
\hspace*{10pt}Submitted to IEEE Trans.\ on IT: CLN 04-326, June 9, 2004; revised Dec.\ 2004, Feb.\ 2005, June 2005. Research partially supported by the NSF and the AFOSR.}}
\date{}
\markboth{
}
{Georgiou: Relative Entropy \& multi-dimensional moment problems}
\maketitle

\begin{abstract}
Entropy-like functionals on operator algebras have been studied since the pioneering work of von Neumann, Umegaki,  Lindblad, and Lieb. The most well-known are the von Neumann entropy
$\mI(\rho):=-\trace (\rho\log \rho)$ and a generalization of the Kullback-Leibler distance
$\mS(\rho ||\sigma ):=\trace (\rho \log \rho - \rho \log \sigma)$, refered to as quantum relative entropy and used to quantify distance between states of a quantum system. The purpose of this paper is to explore $\mI$ and $\mS$ as regularizing functionals in seeking solutions to multi-variable and multi-dimensional moment problems. It will be shown that extrema can be effectively constructed via a suitable homotopy. The homotopy approach leads naturally to a further generalization and a description of all the solutions to such moment problems. This is accomplished by a renormalization of a Riemannian metric induced by entropy functionals.
As an application we discuss the inverse problem of describing power spectra which are consistent with second-order statistics, which has been the main motivation behind the present work.
\end{abstract}

\begin{keywords}Moment problem, spectral analysis, covariance matching, multi-variable, multi-dimensional, quantum entropy.
\end{keywords}

\section{\bf Introduction}
\PARstart{T}{he} quantum relative entropy (Umegaki
\cite{umegaki})
\[\mS(\rho\,\|\sigma):=\trace (\rho \log \rho - \rho \log \sigma)\]
where $\rho,\sigma$ are positive Hermitian matrices (or operators) with trace equal to one, generalizes the
Kullback-Leibler relative entropy \cite{kullbackleibler}, just as the
von Neumann entropy \cite{vonNeumann}
\[
\mI(\rho):=-\trace( \rho \log \rho)
\]
generalizes the classical Shannon entropy. They both inherit a rather rich structure from their scalar counterparts and in particular, $\mS(\cdot \|\cdot)$ is jointly convex in its arguments as shown by Lieb \cite{Lieb} in 1973, whereas $\mI(\cdot)$ is concave. The relative entropy originates in the quest to quantify the difficulty in discriminating between probability distributions and can be thought as a distance between such. 
Its matricial counterpart $\mS$ can similarly be used to quantify distances between positive matrices.

Entropy and relative entropy have played a central  r\^{o}le in thermodynamics in enumerating states consistent with data and, thereby, used to identify ``the most likely'' ones among all possible alternatives. The measurement of a physical property in a classical setting is modeled via {\em ensemble averaging} (e.g., see \cite[Chapter 3]{kittel}) 
\[
r=\sum_k g(k) \rho(k)
\]
where $k$ runs over all micro-states corresponding to a scalar value $g(k)$. Each micro-state occurs with probability $\rho(k)$ and $r$ is a moment of the underlying probability distribution.
Similarly, quantum measurement is also modeled by averaging (as originally idealized by von Neumann, see e.g., \cite[Chapter 5]{sudbery}, \cite[page 183]{halliwell}):
\[
\rho_{\rm after}=\sum_k G(k) \rho_{\rm before} G(k)^*
\]
where the $\rho$'s represent density matrices (positive Hermitian with trace one), the $G$'s represent products of projection operators, and ``$^*$'' denotes ``conjugate-transpose''.
Similar expressions arise for the density operator when restricted to a subsystem (partial trace \cite[page 185]{sudbery}) and  also when measuring ``non-selfadjoint observables'' (e.g., \cite{yen}). If the underlying Hilbert space is infinite dimensional then the measurement process can be modeled with a continuous analogue of the above where the summation is replaced by an integral (e.g., see \cite{braginsky}).
These are instances of moment problems.
More generally we may consider
\begin{equation}\label{eq:discretemoment}
R=\sum_k G_\rleft(k) \rho(k) G_\rright(k)
\end{equation}
where $\rho(k)$ are Hermitian positive matrices
as well as its ``continuous'' counterpart
\begin{equation}\label{generalmomentproblem}
R=\int_\cS G_\rleft(\theta) \rho(\theta)  G_\rright(\theta) d\theta
\end{equation}
where $\rho(\theta)$ represents a Hermitian-valued positive (density) function on a support set $\cS\subseteq \mR^k$ ($k>1$) and $G_\rleft,G_\rright$ are matrix-valued functions on $\cS$. If the underlying distribution is not absolutely continuous then we write
$R=\int_\cS G_\rleft(\theta) d\mu(\theta)  G_\rright(\theta)$ instead,
where $d\mu$ is such a positive Hermitian-valued measure.

The moment problem  (\ref{eq:discretemoment}-\ref{generalmomentproblem}) is typified by spectral analysis based on second-order statistics, especially in the  context of sensor arrays and of polarimetric radar.
The echo at different polarizations and/or at different wavelengths is being sampled at a variety of sensor locations.
It is usually the case that these samples are not independent and that the echo at different frequencies, polarizations etc., affects each sensor by a different amount.
Attributes of the scattering field (e.g., reflectivity at different wavelengths and polarization) and the relative position of the array elements with respect to the scattering field are responsible for the variations in the vectorial echo.
The vector of attributes can be thought of as a vectorial input $u(\theta)$ to the array while the relative position and characteristics of its elements specify a $n_\rleft\times m$ transfer function matrix
\[
G_\rleft(\theta)=\left[\begin{matrix} g_{1,\rleft}\\\vdots\\g_{n_\rleft,\rleft}\end{matrix}\right]
\]
to the $n_\rleft$ sensor outputs.
If the attributes $u(\theta)$ are modeled as a zero-mean vectorial stochastic process, independent over frequencies, then
\[
y_\rleft=\int_\cS G_\rleft(\theta)du(\theta)
\]
represents the vectorial output process. Similarly, if
\[
G_\rright(\theta)=\left[\begin{matrix}g_{1,\rright},&\ldots&g_{n_\rright,\rright}\end{matrix}\right]
\]
is the $m\times n_\rright$ complex conjugate transpose of the transfer matrix corresponding to a second group of sensors,
and if
\[
y_\rright = \int_\cS G_\rright(\theta)^* du(\theta),
\]
designate the corresponding vector of $n_\rright$ outputs, then
the $n_\rleft\times n_\rright$ correlation matrix \[R=E\{y_\rleft y_\rright^*\}\]
gives rise to the matricial moment constraint on the spectral distribution of $u$ given in (\ref{generalmomentproblem}). On the other hand (\ref{eq:discretemoment}) can be interpreted when the power spectrum is discrete.
A power density which matches the correlation samples aims at giving clues about the makeup of the scattering field.

We address the moment problem in the above generality and provide a way to answer the following:
\begin{itemize}{\sf
\item[(i)] does there exist a density function satisfying (\ref{eq:discretemoment}-\ref{generalmomentproblem})?
\item[(ii)] if yes, describe all density functions consistent with (\ref{eq:discretemoment}-\ref{generalmomentproblem}).}
\end{itemize}
\vspace*{.02in}

In essence, the above questions go back to 1980 when
Dickinson \cite{dickinson} raised the issue of consistency of two-dimensional Markovian estimates.
Consistency of {\em scalar} distributions was taken up in the work of Lang and McClellan \cite{LangMcClellan} and Lewis \cite{lewis}.
Both references used entropy functionals and suggested computational solutions to the first question when dealing with scalar distributions.
The present work follows up in the footsteps of these as well as, of a rather extensive literature on
inverse problems \cite{csiszar1,csiszar2,Woods,LevineTribus,miller,dacunha,gamboa,besnerais,junk} having roots in the early days of statistical mechanics. The key idea has been to seek extrema of entropy functionals---existence would guarantee solvability of the moment problem.
The idea of using ``weighted'' entropy functionals to parametrize solutions originates in
Byrnes, Gusev and Lindquist \cite{BGuL}. It was followed up in \cite{BGL1,BGL2} and in \cite{BLkimura,GL} where it was reformulated using the Kullback-Leibler distance between sought solutions and positive ``priors.'' Exploring the connection with the Kullback-Leibler distance, \cite{GL,BLkimura,megretski} and more resently \cite{ac_may2004}, studied the {\em scalar} moment problem in various levels of generality.

Classical moment problems \cite{Akhiezer,KreinNudelman} are closely related to analytic interpolation ones, and as such, have been studied in great generality,
including their matrix-valued counterpart (see e.g., \cite{RR}). However, analytic interpolation applies only when the integration kernels possess a very particular shift-structure similar to that of a Fourier vector (see e.g., \cite{Arov,DD} and also \cite{acmatrix1,acmatrix2}), and is of limited use in the generality sought herein.
On the other hand, literature on interpolation with a ``complexity constraint'' is relevant since it departs from the groove of the classical theory.
Two works are especially relevant, \cite{thesis} and more recently \cite{Blomqvist}. In  \cite{thesis} a homotopy was suggested in the context of the matricial trigonometric moment problem. Then \cite{Blomqvist} used optimization of an entropy functional in the context of matricial Nevanlinna-Pick interpolation. Neither applies in the generality sought herein, yet, below, we build on both of these general directions.

The present work follows up along the lines of \cite{GL} where it was suggested that the quantum relative entropy may be used in the multi-variable case. 
In fact, we carry out the plan suggested in \cite{GL} for multi-variable as well as multi-dimensional distributions and we
develop a computational approach analogous to one presented in \cite{ac_may2004} for scalar distributions.

In view of a rather rich literature on quantum entropies, a comment is in order as to other possible connections to the present work.
Besides the Umegaki-von Neumann entropy $\mS(\cdot \|\cdot)$ studied in this paper, there is a plethora of alternatives due to a dichotomy between matricial and scalar distributions \cite{petzsudar,Ohya}. In particular Araki's theory \cite{araki,ruskai} helped characterize a family of ``quasi-entropies,'' contractive under stochastic maps. References \cite{petz1,petz2,ruskai} in particular explore the Riemannian geometry they induce on density matrices. It is an interesting question as to which among this ``garden of entropies,''  besides the Umegaki-von Neumann one, allows a convenient representation of solutions for general moment problems. The approach we have taken leads us to work mostly with an induced metric (a Jacobian related to the Hessian of $\mS(\cdot \|\cdot)$). A suitable normalization then recovers any solution of the moment problem as a corresponding extremal. It is not known whether, the ``weighted metrics'' e.g., $\nabla h_W$ in Section \ref{completesolutions}, are metrics induced by a quasi-entropy in the language of Petz \cite{petz1}.

Finally, it is interesting to point out that, a counterpart for discrete distributions relates to the theory of analytic centers in semi-definite programming \cite{sboyd,NN}. In fact, a key construction in this paper---a homotopy for the numerical computation of solutions, is analogous to tracing paths of analytic centers in interior point methods.

In Section \ref{section:motivating} with discuss three motivating examples. Section \ref{section:matricial}
develops the geometry of matricial cones and the significance of relative entropy as a barrier functional. In order to simplify the exposition, we first deal with cones of constant matrices and then with those of matrix-valued density functions. Except for technical differences, the two run in parallel. Finally, Section \ref{completesolutions} discusses the parametrization of solutions to the moment problem,
followed by concluding remarks (Section \ref{conclusion}) and an appendix (Section \ref{appendix})
with useful facts on matricial calculus.

\section{\bf Motivating Examples}\label{section:motivating}

In this section we present four examples that motivate our study. For simplicity, the first two are developed in the context of scalar distributions, the third is intrinsically multi-variable. The fourth example pertains to the connection of the moment problem with analytic interpolation. We return to the first and fourth example again later on in the paper.

\subsection{\bf Non-equispaced arrays}\label{example1}
\begin{figure}[htb]\begin{center}
\psfrag{phi}{$\phi$}
\psfrag{E1}{$E_0$}
\psfrag{E2}{$E_1$}
\psfrag{E3}{$E_2$}
\includegraphics[totalheight=8cm]{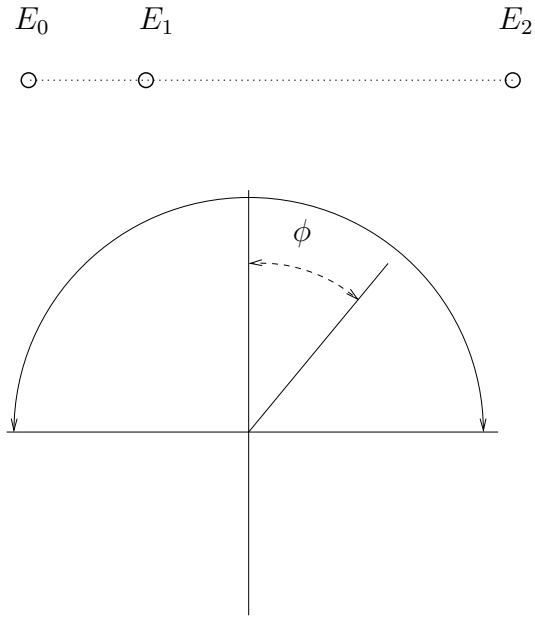}
\caption{Non-equispaced sensor array}\label{fig:sensarray1}\end{center} \end{figure}

Consider an array of sensors with three elements, linearly spaced at distances $1$ and $\sqrt{2}$ wavelengths from one another, and assume that (monochromatic) planar waves, originating from afar, impinge upon the array. This is exemplified in Figure \ref{fig:sensarray1}.

Assuming that the sensors are sensitive to disturbances originating over one side of the array,
with sensitivity independent of direction,
the signal at the $\ell$th sensor is typically represented as a superposition
\[
u_\ell(t)=\int_0^\pi A(\theta) e^{j(\omega t - p x_\ell\cos(\theta) + \phi(\theta))}d\theta,
\]
of waves arising from all spatial directions $\theta\in[0,\pi]$,
where $\omega$ is as usual the angular time-frequency (as opposed to ``spatial''),
$x_\ell$ the distance between the $\ell$th and the $0$th sensor, $p$ the wavenumber, and $A(\theta)d\theta$ the amplitude and $\phi(\theta)$ a random phase of the $\theta$-component. Typically, $\phi(\theta)$ for various values of $\theta$ are uncorrelated.
The term $ p x_\ell\cos(\theta)$ in the exponent
accounts for the phase difference between reception
at different sensors. For simplicity we assume that $p =1$ in appropriate units.
Correlating the sensor outputs we obtain
\[R_k = E\{ u_{\ell_1}\bar{u}_{\ell_2}\}:=  \int_0^\pi e^{-jk \cos(\theta)} f(\theta) d\theta\]
where $f(\theta)=|A(\theta)|^2$ represents {\em power density}, and $k=\ell_1-\ell_2$ with $\ell_1\geq \ell_2$ and belonging to $\{0,1,\sqrt{2}+1\}$. Thus,
\begin{equation}\label{indexingset}
k\in\cI:=\{0,1,\sqrt{2}, \sqrt{2}+1\}.
\end{equation}
The only significance of our selection of distances between sensors, that gave rise to the indexing set (\ref{indexingset}), is to underscore that there is no algebraic dependence between the elements of the so-called {\em array manifold}
\[ G(\theta):=\left[\begin{matrix}1 & e^{-j\tau} & e^{-j\sqrt{2}\tau} & e^{-j(\sqrt{2}+1)\tau}\end{matrix}\right]^\prime,
\]
(where ``$[\cdot ]^\prime$'' denotes the transpose, $G$ is thought of as a column vector, and $\tau=\cos(\theta)\in[-1,1]$).

Given a set of values $R_k$ for $k\in\cI$, it is often important to determine whether they are indeed the moments of a power density $f(\theta)$, and  if so to characterize all consistent power spectra.
The case of  arrays with equispaced elements is very special
and answers to such questions relate to the non-negativity of a Toeplitz matrix formed out of the $R_k$'s. In the present situation nonnegativity of
\[
\int_{-1}^1 \left[\begin{matrix}1 \\ e^{-j\tau} \\ e^{-j\sqrt{2}\tau}\end{matrix}\right]\frac{f(\cos^{-1}(\tau))}{\sqrt{1-\tau^2}}\left[\begin{matrix}1 & e^{j\tau} & e^{j\sqrt{2}\tau}\end{matrix}\right]d\tau
\]
which, in  the obvious indexing turns out to be
\begin{equation}\label{skew_Toeplitz} \left[\begin{matrix}R_0& R_1 & R_{\sqrt{2}+1}\\
\bar{R}_1 & R_0 & R_{\sqrt{2}}\\
\bar{R}_{\sqrt{2}+1} & \bar{R}_{\sqrt{2}}& R_0\end{matrix}\right],
\end{equation}
is  only a necessary condition. The fact that it is {\em not} sufficient (see e.g., \cite[page 786]{sesha1}) motivated the present study.



\subsection{\bf Two-dimensional distributions}\label{example2}

The subject matter of multi-dimensional distributions received considerable attention in the 1970's and early 1980's (e.g., see \cite{marzetta}).
However, despite the rich theory of analytic interpolation and orthogonal polynomials in more than one variable, etc.\ (e.g., see \cite{jackson}, \cite{rudin}), as pointed out by Brad Dickinson \cite{dickinson}, the analog of the questions raised earlier never received a definitive answer.

For instance, consider the simplest possible planar grid
\[
\cI:=\{k,\ell\in\mZ\;:\; 0\leq k\leq n,\; 0\leq \ell \leq n\}
\]
which is both regular and square. In fact, we may even assume that $n=2$. Then,
consider the case where (monochromatic as before) waves impinge upon sensors at the grid nodes,
the sensors are sensitive to reception on one side of their plane, the waves are planar, originating from all possible directions in the sky and uncorrelated over different directions, and that correlations at the sensor outputs are taken. The array manifold in this case becomes (after suitable normalization assumptions) the $3\times 3$ array
\[
G(\theta,\phi)=\left[ e^{j (k\theta+\ell \phi)}\right]_{k,\ell=0}^{k,\ell=2},
\]
with $(\theta,\phi)\in[0,\pi]\times[0,\pi]$.
This is quite standard (e.g., \cite{johnson}, \cite{vantrees}, \cite{Haykin_array}).
Now, if $\rho(\theta,\phi)$ denotes a positive scalar power density for the waves originating from (normalized) Euler angles $\theta,\phi$ then, correlation between the sensor outputs gives us a $3\times 3$ matrix of covariance samples
\[
\left[ R_{k,\ell}:= \int_0^\pi \int_0^\pi e^{j (k\theta+\ell \phi)}\rho(\theta,\phi)d\theta d\phi \right]_{k,\ell=0}^{k,\ell=2}.
\]
The same two questions that we raised earlier are again relevant:
Given a $3\times 3$ array, how can we tell that it originates as suggested above,
and how can we characterize {\em all} power densities that are be consistent with the covariance samples $R_{k,\ell}$?

Lang and McClellan \cite{LangMcClellan} considered ``maximum entropy spectra'' and suggested evaluating those on a discrete grid\footnote{Lang and McClellan \cite{LangMcClellan} were the first to point out that existence cannot be guaranteed for the usual ``entropy rate'' functional $\int_\cS \log(\rho)d\theta$ unless the dimension of $\cS$ is one, cf.\ Theorem \ref{thm1:matrix} below.}. Lewis \cite{lewis} considered a framework which is applicable for answering the existence question, while
\cite{ac_may2004} discussed parametrization of solutions as well. The present contributions allows addressing the most general situation  where the sensor array elements receive vectorial echoes and hence, $R_{k,\ell}$ as well as $\rho(\theta,\phi)$ are matrices.

\subsection{\bf Quantum measurements}\label{example3}

We temporarily adopt the language of quantum mechanics as in, e.g., \cite{NielsenChuang}, and explain how this relates to the linear algebraic framework in the introduction. We then discuss a basic academic paradigm which exemplifies the setting of the matricial moment problem.

Let $\rho^{AB}$ denote the density matrix of a quantum system composed of two subsystems $A$ and $B$. Each subsystem can be in two states $|0\rangle$ and $|1\rangle$, respectively. These can be thought of as the vectors
\[
\left(\begin{matrix} 1\\ 0\end{matrix}\right)\mbox{ and } \left(\begin{matrix} 0\\ 1\end{matrix}\right) \mbox{ in }\mC^2.
\]
Then, the states of the combined system, $| i\rangle\otimes | j\rangle$ where $i,j\in\{0,1\}$, can be represented by a vector in $\mC^4$. For instance, $|0\rangle\otimes|1\rangle$ corresponds to
\[ \left(\begin{matrix} 1\\ 0\end{matrix}\right)\otimes \left(\begin{matrix} 0\\ 1\end{matrix}\right) =\left(\begin{matrix}0\\ 1\\ 0\\ 0\end{matrix}\right)
\]
where in the last equation $\otimes$ is the Kronecker tensor product. Accordingly, $\rho^{AB}$ is a Hermitian $4\times 4$ trace one matrix formed out of sums of Kronecker products of $2\times 2$ density matrices of the two subsystem. Now, if $\rho^A$ represents the density matrix of subsystem $A$, then this can be obtained from $\rho^{AB}$ via taking the partial trace with respect to subsystem $B$ (e.g., see \cite[pages 106-107]{NielsenChuang}).
An important example is that of the Bell state with
\[
\rho^{AB}=\left(\frac{|01\rangle +|11\rangle}{\sqrt{2}}\right)\left(\frac{\langle 00|+\langle 11|}{\sqrt{2}}\right) =
\frac12\left(\begin{matrix} 1 & 0 & 0 & 1\\
                                              0 & 0 & 0 & 0\\
                                              0 & 0 & 0 & 0\\
                                              1 & 0 & 0 & 1\end{matrix}\right).
\]
The partial trace then gives (see \cite[pages 106-107]{NielsenChuang})
\[
\rho^A = \trace_B(\rho^{AB})= \frac12 I_2
\]
where $I_2$ is the $2\times 2$ identity matrix. If $\rho^{AB}$ is represented as $4\times 4$ matrix as above, then
\[
\rho^{A}= \sum_{k=1}^2 G_k \rho^{AB} G_k^*
\]
where
\[
G_1= I_2\otimes \left(\begin{matrix} 1& 0\end{matrix}\right) =\left(\begin{matrix} 1 & 0 & 0 & 1\\
                                                                                         0 & 0 & 1 & 0\end{matrix}\right),
 \]
 and similarly, $G_2=I_2\otimes \left(\begin{matrix} 0 & 1\end{matrix}\right)$.
The reverse task of characterizing all possible $\rho^{AB}$ based on a known $\rho^A$ is a (discrete) moment problem.                                               

\subsection{\bf State-statistics and analytic interpolation with degree constraint}\label{example4}
Consider the linear discrete-time state equations 
\begin{eqnarray}
x_k&=&A x_{k-1} + B u_k,                                  \label{i2s}
\mbox{ for $k\in{\mathbb Z}$},
\end{eqnarray}
where $x_k\in\mC^n,$ $u_k\in\mC^m$,
$A\in\mC^{n\times n},$ $B\in\mC^{n\times m}$, $(A,B)$ is a controllable pair, and
the eigenvalues of $A$ lie in the open unit disk of the complex plane.
Let $\{u_k:k\in{\mathbb Z}\}$ be a zero-mean stationary stochastic process
with power spectrum 
the non-negative matrix-valued measure $d\mu(\theta)$ on $\theta\in(-\pi,\pi]$.
Then, {\em under stationarity conditions,} the state covariance
\[
R:=E\{x_kx_k^*\}
\]
can be expressed in the form of the integral (cf.\ \cite[Ch.\ 6]{masani})
\begin{equation}\label{Sigma}
R=\int_{-\pi}^{\pi}
\left( G(e^{j\theta}) \frac{d\mu(\theta)}{2\pi} G(e^{j\theta})^* \right)
\end{equation}
where
\[
G(z):=(I-z A)^{-1}B
\]
is the transfer function of system (\ref{i2s}).
Note that we use $z$ to denote the transform of the delay operator and therefore
$G(z)$ is analytic in the unit disc of the complex plane.

It turns out that state covariances $R$ are characterized by the following two equivalent conditions (see \cite{acmatrix1,acmatrix2})
\begin{equation}\label{rankcondition}
{\rm rank}\left[\begin{array}{cc}R-AR A^*&B\\B^*&0\end{array}\right]
=
2m
\end{equation}
and,
\begin{equation}\label{sylvester}
R-AR A^*=BH+H^*B^* \mbox{ for some }H\in\mC^{m\times n}.
\end{equation}
Then, power spectral measures consistent with (\ref{Sigma}) are in correspondence with matrix valued functions $F(z)$ on the unit circle ${\mathbb{D}}:=\{z\in{\mathbb{C}} : |z|<1\}$ which have nonnegative real part via the Herglotz representation
\begin{equation}\label{correspondingpr}
F(z)=
\int_{-\pi}^{\pi}\left(\frac{1+z e^{j\theta}}{1-z e^{j\theta}}\right)\frac{d\mu(\theta)}{2\pi} + jc,
\end{equation}
with $jc$ an arbitrary skew-Hermitian constant.
The measure $d\mu$ can be recovered as the weak* limit of
the real part of $F(z)$ as $z$ tends to the boundary, i.e.,
\begin{equation}\label{realpart}
d\mu(\theta)\sim \lim_{r\nearrow 1}{\Re}(F(re^{j\theta})).
\end{equation}
The class of nonnegative real matrix valued functions $F$ giving rise to admissible power spectral measures are also characterized by the interpolation  condition (\cite{acmatrix1})
\begin{equation}\label{PRNehari}
F(z)=H(I-z A)^{-1}B + Q(z) V(z)
\end{equation}
where $Q$ is a matrix function analytic in $\mathbb D$,
\begin{equation}\label{V}
V(z):=D+z C(I-z A)^{-1}B
\end{equation}
and  $C\in\mC^{m\times n}$, $D\in\mC^{m\times m}$ are selected so that
$V$ is inner, i.e., $V(\xi)^*V(\xi)=V(\xi)V(\xi)^*
=I$ for all $|\xi|=1$.

The data $A,B,H$ and $V(z)$ in equation (\ref{PRNehari}) specify an analytic interpolation problem. Positive-real solutions to (\ref{PRNehari}) can be given via (\ref{correspondingpr}) and solutions to the moment problem (\ref{Sigma}). The characterization of solutions to (\ref{Sigma}) given in Theorem \ref{thm1:matrixW} allows
a non-classical characterization of solutions to (\ref{PRNehari}) and in particular a characterization of solutions of McMillan degree less than or equal to the dimension of (\ref{i2s}) (see Section \ref{example4_continued}).

\section{\bf Matricial distributions and their moments}\label{section:matricial}

The moment conditions (\ref{eq:discretemoment}-\ref{generalmomentproblem}) are  linear constraints on densities $\rho_k$ ($k=1,2,\ldots$) and $\rho(\theta)$ ($\theta\in\cS$), respectively. Density functions, whether discrete or continuous, are non-negative, or non-negative definite in the matricial case, for each value of their indexing set. Thus they have the structure of a cone.
Entropy functionals on the other hand represent natural barriers on such positive cones and can be used to identify, and even parametrize, density functions which are consistent with given moment conditions. We begin by explaining the geometry of the moment problem for constant density matrices
and the relevance of entropy functionals in obtaining solutions
as their respective extrema. Both, the geometry of cones of matricial densities functions
as well as the r\^{o}le of entropy functionals
is quite similar and is taken up in Section \ref{sec:matricialdistributions}.

\subsection{\bf Relative entropy and the geometry of matricial cones}\label{sec:constantdensities}

We begin by focusing on constraints
\[
R=\sum_k G_\rleft(k) \rho\, G_\rright(k)
\]
where $\rho$ is not indexed. The general case is quite similar.

We use the notation
\newcommand{\mM}{{\mathbb M}}\newcommand{\dual}{{\rm dual}}
\begin{eqnarray*}
\mM&:=&\{M\in\mC^{m\times m}\;:\;M=M^*\},\\
\fM&:=&\{M\in \mM \mbox{ and } M\geq 0\},\\
\fM_+&:=&\{M\in \mM \mbox{ and } M> 0\}
\end{eqnarray*}
to denote the space of Hermitian matrices and the cones of non-negative and positive definite ones, respectively.
The space $\mM$ is endowed with a natural
inner product
\[ \langle M_1,M_2\rangle := \trace(M_1^* M_2)=\trace(M_1M_2)
\]
as a linear space over $\mR$.
Clearly, both, $\fM$ and $\fM_+$ are convex cones.
Since non-negativity of $\langle M_1,M\rangle$ for all $M_1\in\fM$ implies that $M\in\fM$, it follows that $\fM$ is self-dual\footnote{In general, the dual cone $\fM^\dual$ is the set of elements forming an  ``acute angle'' with all elements of the original cone, i.e., $\{M\;:\;\langle M,M_1\rangle\geq 0, \;\forall M_1\in\fM\}$ (see \cite{KreinNudelman}).}. It can also be seen that
$\fM_+$ is the interior of $\fM$.

The linear operator
\[
L : \mM \to \fR : \rho \mapsto R=\sum_k G_\rleft(k) \rho\, G_\rright(k)
\]
where $\fR\subseteq\mC^{n_\rleft \times n_\rright}$ denotes the range of $L$,
maps $\fM$ onto the cone of admissible moments $\cK= L(\fM)\subseteq \fR$. Here, and throughout, $G_\rleft,G_\rright$ are matrices of dimension $n_\rleft\times m$ and $m\times n_\rright$, respectively. A further assumption that is often needed is that the null space of $L$ does not intersect $\fM$, i.e.,
\begin{equation}\label{nonintersect}
{\rm null}(L)\cap \fM=\{0\}.
\end{equation}
The interior of $\cK$ is ${\rm int}(\cK)=L(\fM_+)$ and, 
given $R$, the moment problem requires testing whether
$R\in\cK$ and if so, characterizing all $\rho\in\fM$ such that $R=L(\rho)$.

Geometry in the range space $\fR$ is based on
\begin{equation}\label{innerproduct}
\langle \lambda, R\rangle := \Real \left(\trace(\lambda^* R)\right), \mbox{ for }\lambda, R\in\fR.
\end{equation}
Then the adjoint transformation of $L$ is
\[
L^* : \fR \to \mM : \lambda \mapsto \rho=\left(\sum_k G_\rleft(k)^* \lambda\, G_\rright(k)^*\right)_{\rm Herm}
\]
where
$(M)_{\rm Herm}:=\frac12 \left( M+M^*\right)$ is the ``Hermitian part''.
The dual cone of $\cK$,
\[
\cK^\dual:= \{ \lambda \in\fR\;:\; \langle \lambda,R\rangle \geq 0, \;\forall R\in\cK\},
\]
is naturally related to the cone $\fM\subset\mM$. In fact, using $\langle \lambda,L(\rho)\rangle=\langle L^*(\lambda),\rho\rangle$ it follows easily that
\[
\cK^\dual=\{\lambda \in\fR\;:\; L^*(\lambda) \in \fM\}.
\]
The interior of the dual cone
\[
\cK_+^\dual:=\;{\rm int}(\cK^\dual):=\left\{\lambda\;:\; \langle \lambda, R\rangle >0,\; \forall R\in\cK-\{0\}\right\}
\]
corresponds to $\fM_+$ as is easily seen to satisfy
\[
\cK_+^\dual=\{\lambda\;:\; L^*(\lambda)\in\fM_+\}.
\]
Finally, (\ref{nonintersect}) can be seen to be equivalent to $\cK_+^\dual\neq\emptyset$.

\subsubsection{\bf Minimizers of $\mS(I\|\rho)$:}\label{sec:minimizersofrho}
We are interested in minimizers of (the negative entropy)
\[
\mS(I\|\rho)=-\trace(\log(\rho))
\]
on $\fM_+$ subject to $R=L(\rho)$. Here and throughout, ``$I$'' denotes the identity matrix of size determined from the context.
When such a minimizer exists at an interior point of $\fM_{R,+}$, stationarity conditions for the entropy functional
dictate an explicit form for the minimizer (which, is unique due to the convexity of $-\trace(\log(\rho))$).

The Lagrangian of the problem is
\newcommand{\cL}{{\cal L}}
\[
\cL(\lambda,\rho):= \trace(-\log(\rho))-\langle\lambda,R-L(\rho)\rangle.
\]
Using the expression for the derivative of the logarithm given in Proposition \ref{diffoflog} of the appendix,
the (Gateaux) derivative of $\cL$ in the direction $\delta\in\mM$ becomes
\begin{eqnarray*}
d \cL(\lambda,\rho \,;\, \delta)&:=& \trace(-M_\rho^{-1} \delta)+ \langle \lambda,L(\delta)\rangle\\
&=& \trace(-\rho^{-1}\delta) + \langle L^*(\lambda),\delta\rangle.
\end{eqnarray*}
In the above derivation, the  ``trace'' is what allows replacing the ``non-commutative division operator'' $M_\rho^{-1}$ (cf.\ (\ref{eigenexpansion})) with multiplication by $\rho^{-1}$.
The stationarity condition $d \cL(\lambda,\rho\,;\, \delta)\equiv 0$ then gives
\begin{eqnarray}\label{logoptimal}
\rho&=& \left( L^*(\lambda) \right)^{-1}. \label{formofoptimum}
\end{eqnarray}
Thus, a necessary condition is that there exist $\lambda\in\cK^\dual$ such that
$L^*(\lambda)$ is strictly positive, i.e., that $\cK^\dual_+$ is nonempty.
It turns out that if $R\in{\rm int}(\cK)$ then this condition is also sufficient.

\begin{thm}\label{thm1a}{\sf Assume that $R\in{\rm int}(\cK)$. Then the entropy functional $\mS(I\|\rho)$
has a minimum in $\fM_{R,+}$, which is also unique, if and only if $\cK_+^\dual$ is nonempty.
}\end{thm} 

\begin{proof} Necessity is obvious and was argued above. Uniqueness follows from the matrix convexity of $\log(\cdot)$ (which follows e.g., from a positive Hessian (\ref{loghessian}) given in the appendix). Sufficiency requires that there exists $\lambda_1\in\cK^\dual_+$ such that $R=L(\left( L^*(\lambda_1) \right)^{-1})$. Then
$
\rho_1=\left( L^*(\lambda_1) \right)^{-1}
$
satisfies both the stationarity conditions and the contstraints.
We show this via a continuity argument which we will adopt again later on to more general cases.

Consider the mapping
\begin{eqnarray*}
h&:& \cK_+^*\to {\rm int}(\cK)\subset \fR\\
&:& \lambda \mapsto L(\left( L^*(\lambda) \right)^{-1}).
\end{eqnarray*}
Its Jacobian
\begin{eqnarray*}
\nabla h |_\lambda &:&\fR\to\fR \;:\; \delta \mapsto L( L^*(\lambda)^{-1}  L^*(\delta) L^*(\lambda)^{-1}).
\end{eqnarray*}
is Hermitian since
\begin{eqnarray*}
\langle \delta_1, \nabla h (\delta)\rangle &=& \langle L^*(\delta_1),L^*(\lambda)^{-1}L^*(\delta)L^*(\lambda)^{-1}\rangle\\
&=&  \langle L^*(\lambda)^{-1}L^*(\delta_1)L^*(\lambda)^{-1},L^*(\delta)\rangle\\
&=&\langle \nabla h (\delta_1),\delta\rangle.
\end{eqnarray*}
The map $h$ is a local diffeomorphism and $\nabla h$ can be used to relate locally a smooth path in the space of $R$'s to 
one in the  space of $\lambda$'s (of course, both spaces being the same space $\fR$).

Choose a $\lambda_0\in\cK_+^*$, let
\[
\rho_0=\left( L^*(\lambda_0) \right)^{-1},
\]
and let $R_0=L(\rho_0)$.
Consider the interval path $R_\tau=(1-\tau)R_0+\tau R$, for $\tau\in[0,1]$. Since $R\in{\rm int}(\cK)$, so is the whole path $R_\tau$ ($\tau\in[0,1]$). We claim that for all $\tau\in[0,1]$ there exists
$\lambda_\tau\in\cK^\dual_+$ such that $R_\tau=L(L^*(\lambda_\tau)^{-1})$. It is clear that this holds locally
and that $\lambda_\tau$ satisfies
\begin{equation}\label{firstoccurence}
\frac{d}{d\tau} \lambda_\tau =  (\nabla h|_{\lambda_\tau})^{-1}(R-R_0),
\end{equation}
since $\frac{d}{d\tau}R_\tau=R-R_0$.
The starting point is $\lambda_0$ and (\ref{firstoccurence}) can be integrated
over a maximal interval $[0,\epsilon)$ for which $\lambda_\tau\in\cK^\dual_+$. Throughout
\[
R_\tau=L(L^*(\lambda_\tau)^{-1}).
\]
If $\epsilon>1$, this proves our claim. If $\epsilon\leq 1$, then either $\|\lambda_\tau\|\to\infty$ as $\tau\to\epsilon$, or the $\lambda_\tau$'s have a limit
point $\lambda_\epsilon$ on the boundary of $\cK^\dual_+$, i.e., such that $L^*(\lambda_\epsilon)$ is singular.
Below we argue that neither is possible, which then shows that $\epsilon>1$ and completes the proof.

We first show that $\lambda_\tau$ remains bounded. Assume to the contrary, i.e., assume that $\|\lambda_\tau\|$ grows unbounded as $\tau\to\epsilon$,
and let $\ell_\tau:=\lambda_\tau/ \|\lambda_\tau \|$ (where $\|\lambda\|=\sqrt{\langle\lambda,\lambda\rangle}$ as usual).
Since $R_\epsilon\in{\rm int}(\cK)$, it holds that $R_\epsilon=L(\rho)$ with $\rho>0$, and $\langle \ell_\tau,R\rangle$ is bounded away from zero for elements $\ell_\tau\in\cK^\dual$ of unit norm.
However, because
\[
\langle \lambda_\tau,R_\tau\rangle = \langle L^*(\lambda_\tau), L^*(\lambda_\tau)^{-1}\rangle = \trace(I),
\]
it follows that $\langle \ell_\tau,R_\tau\rangle=\trace(I)/\|\lambda_\tau\|$ has zero as a limit point when $\tau\in[0,\epsilon)$; hence, so does $\langle \ell_\tau,R_\epsilon\rangle$. But this is a contradiction, hence $\lambda_\tau$ remains bounded.

We finally show that $L^*(\lambda_\tau)^{-1}$  and $\nabla h$ along with its inverse remain bounded.
Consider the quadratic form
\begin{equation}\label{quadraticform}
\langle \delta,\nabla h_{\lambda_\tau}(\delta)\rangle = \|L^*(\lambda_\tau)^{-1/2}L^*(\delta)L^*(\lambda_\tau)^{-1/2}\|^2,
\end{equation}
for $\delta\in\fR$.
Since $\lambda_\tau$ (and hence $L^*(\lambda_\tau)$) remains bounded, the quadratic form is bounded away from  zero
when $\tau\in[0,\epsilon)$.  Hence, $\nabla h|_{\lambda_\tau}^{-1}$ is uniformly bounded on $[0,\epsilon)$.
On the other hand, because of (\ref{nonintersect}), the minimal angle between
any ray in the cone $\fM$ and ${\rm range}(L^*)$ is bounded away from $\pi/2$.
Hence, $\|R_\tau\|> \alpha \|L^*(\lambda_\tau)^{-1}\|$, for some $\alpha>0$. But $R_\tau$ remains bounded.
We conclude that $L^*(\lambda_\tau)^{-1}$ remains bounded and that $\nabla h$ remains bounded as well.
This completes the proof.
\end{proof}

\subsubsection{\bf Minimizers of $\mS(\rho\|I)$:}\label{expconstant}

We now focus on minimizers of
\[
\mS(\rho\|I)=\trace\left( \rho\log(\rho)\right)
\]
in $\fM_+$, subject to $R=L(\rho)$.
The Lagrangian this time is
\[
\cL(\lambda,\rho):= \trace(\rho\log(\rho))-\langle\lambda,R-L(\rho)\rangle.
\]
Once again, using the expression for the differential of the logarithm
given in Proposition \ref{diffoflog} of the appendix,
the (Gateaux) derivative of $\cL$ in the direction $\delta\in\mM$ becomes
\begin{eqnarray*}
d \cL(\lambda,\rho \,;\, \delta)&:=& \trace(\delta\log(\rho)+\rho M_\rho^{-1}(\delta))+ \langle \lambda,L(\delta)\rangle\\
&=& \trace(\delta\log(\rho)+\delta) + \langle L^*(\lambda),\delta\rangle.
\end{eqnarray*}
The last step follows from
\begin{eqnarray*}
\trace(\rho M_\rho^{-1}(\delta))&=&\trace(\rho \int_0^\infty (\rho+t)^{-1}\delta (\rho+t)^{-1}dt)\\
&&\hspace*{-60pt}=\trace(\int_0^\infty (\rho+t)^{-1}\rho(\rho+t)^{-1}dt\,\delta)=\trace(\delta).
\end{eqnarray*}
The stationarity condition $d \cL(\lambda,\rho\,;\, \delta)\equiv 0$ then gives that
\begin{equation}
\rho=\exp\left(-I- L^*(\lambda)\right)= \frac{1}{e}\exp\left(-L^*(\lambda)\right) \label{formofoptimumexp}
\end{equation}
with $L^*(\lambda)\in\mM$ (and not necessarily in $\fM$ as before).
It turns out that if $R\in{\rm int}(\cK)$ a minimizer can always be found.
It should be noted that (\ref{nonintersect}) is no longer a necessary condition.

\begin{thm}\label{thm2a} {\sf If $R\in{\rm int}(\cK)$, then the entropy functional $\mS(\rho\| I)$ has a minimum in $\fM_{R,+}$ which is unique and of the form
(\ref{formofoptimumexp}).
}\end{thm} 

\begin{proof}
We use a similar continuity argument as before.
Consider the mapping
\begin{eqnarray}\label{kappaone}
\kappa &:& \fR\to {\rm int}(\cK)\subset \fR\\
&:& \lambda \mapsto L\left(\frac{1}{e} \exp(-L^*(\lambda)) \right).\nonumber
\end{eqnarray}
Its Jacobian
\begin{equation}\label{nablakappaone}
\nabla k |_\lambda \;:\;\fR\to\fR \;:\; \delta \mapsto \frac{1}{e} L\left(   M_{\exp(-L^*(\lambda))}(-L^*(\delta))\right),
\end{equation}
with $M_C$ as in (\ref{multiplicationoperator}), is Hermitian and negative definite. This is because
\begin{eqnarray*}
\langle \delta_1, \nabla \kappa|_\lambda (\delta)\rangle &=& -\frac{1}{e} \langle L^*(\delta_1), M_{\exp(-L^*(\lambda))}(L^*(\delta))\rangle\\
&& \hspace*{-60pt}= -\frac{1}{e}  \trace( L^*(\delta_1) \int_0^1e^{-(1-t)L^*(\lambda)} L^*(\delta)e^{-tL^*(\lambda)}dt )\\
&& \hspace*{-60pt}= -\frac{1}{e}  \trace(\int_0^1e^{-tL^*(\lambda)}L^*(\delta_1)e^{-(1-t)L^*(\lambda)}\,dt L^*(\delta))\\
&=&\langle \nabla \kappa|_\lambda (\delta_1),\delta\rangle,
\end{eqnarray*}
while
\[\langle \nabla \delta,\kappa|_\lambda (\delta)\rangle=
 -\frac{1}{e}  \int_0^1  \trace(A^\frac{t}{2} L^*(\delta)  A^{-\frac{t}{2}} A A^{-\frac{t}{2}}L^*(\delta)  A^{\frac{t}{2}}  )     dt
\]
with $A=\exp(-L^*(\lambda))\in\fM_+$. Then $\langle \nabla \delta,\kappa|_\lambda (\delta)\rangle$ is clearly negative unless $\delta=0$.
Thus, the map $\kappa$ is a local diffeomorphism and $\nabla \kappa$ can be used to relate
locally a smooth path in the space of $R$'s to 
one in the  space of $\lambda$'s.

Begin with $\lambda_0$ in $\fR$, $R_0=L(\frac{1}{e}\exp(-L^*(\lambda_0))$, and
$R_\tau=(1-\tau)R_0+\tau R$ for $\tau\in[0,1]$. Since $R_0,R\in{\rm int}(\cK)$ then $R_\tau\in{\rm int}(\cK)$ for all $\tau\in[0,1]$. We claim that
\begin{equation}\label{firstoccurenceexp}
\frac{d}{d\tau} \lambda_\tau =  (\nabla \kappa|_{\lambda_\tau})^{-1}(R-R_0),
\end{equation}
can be integrated over $[0,1]$ with $\lambda_\tau$ staying bounded. Then, by construction,
$\lambda_\tau$ satisfies both
\begin{equation}\label{Rlambdaopen}
R_\tau=L(\frac{1}{e}\exp(-L^*(\lambda_\tau))
\end{equation}
as well as the stationarity conditions for each $\tau$. Hence,
\[
\rho=\frac{1}{e}\exp(-L^*(\lambda_1)
\]
is a minimizer for $\mS(\rho\|I)$ as claimed in the proposition.

Clearly (\ref{firstoccurenceexp}) can be integrated over $[0,\epsilon)$. If $\epsilon>1$, we are done.
Thus we only need to show that $\epsilon \leq 1$ and
\[
\|\lambda_\tau\|\to\infty \mbox{ as }\tau\to\epsilon
\]
lead to a contradition. To this end, we use two facts, first that (\ref{Rlambdaopen}) holds on $[0,\epsilon)$
and then, that $\langle \ell,R_\epsilon\rangle>0$ for all $\ell\in\cK^\dual$, since $R_\epsilon\in{\rm int}(\cK)$.

Define $\ell_\tau=\lambda_\tau/\|\lambda_\tau\|$ and let $\ell_\epsilon$ be a limit point of a convergent subsequence $\ell_{\tau_i}$ with $\tau_i\to \epsilon$ (which exists since the $\ell_\tau$'s are bounded).
We claim that $L^*(\ell_\epsilon)\geq 0$ and singular. To see this first note that, if 
\[
{\rm spectrum}(L^*(\ell_\epsilon))\cap (-\infty,0)\neq \emptyset,
\]
then ${\rm spectrum}(L^*(\ell_{\tau_i}))\cap (-\infty,0)\neq \emptyset$ as well, for sufficiently large $i$.
But if this is so, then
\begin{equation}\label{zeroinf}
L(\frac{1}{e}e^{-L^*(\lambda_{\tau_i})})=L(\frac{1}{e}e^{-L^*(\ell_{\tau_i}) \| \lambda_{\tau_i}\|})
\end{equation}
grows without bound instead of tending to $R_\epsilon$. Therefore,
\[
{\rm spectrum}(L^*(\ell_\epsilon))\subset [0,\infty).
\]
Now, if $L^*(\ell_\epsilon)$ is nonsingular then for all $i$ sufficiently large $L^*(\ell_{\tau_i})$ is nonsingular as well. But then, (\ref{zeroinf}) tends to zero as $i\to\infty$. Hence, $0\in{\rm spectrum}(L^*(\ell_\epsilon))\subset [0,\infty)$.

Now let $U$ be an isometry ($U^*U=I$) whose columns span the range of $L^*(\ell_\epsilon)$
and consider
\begin{eqnarray*}
\langle \ell_\epsilon,R_{\tau_i}\rangle&=&\frac{1}{e}\langle L^*(\ell_\epsilon), e^{-L^*(\lambda_{\tau_i})})=\frac{1}{e}\langle L^*(\ell_\epsilon), e^{-L^*(\ell_{\tau_i}) \|\lambda_{\tau_i}\|})\\
&=&\frac{1}{e}\langle U^*L^*(\ell_\epsilon)U, e^{-U^*L^*(\ell_{\tau_i})U \|\lambda_{\tau_i}\|}).
\end{eqnarray*}
Since 
\[U^*L^*(\ell_{\tau_i})U\to U^*L^*(\ell_\epsilon)U>0
\]
while $\|\lambda_{\tau_i}\|\to\infty$ we conclude that $\langle \ell_\epsilon,R_{\tau_i}\rangle\to 0$
as $\tau\to\infty$. Since, $R_{\tau_i}\to R_\epsilon$ it follows that
$\langle \ell_\epsilon,R_\epsilon\rangle=0$ which contradicts the hypothesis that $R_\epsilon\in{\rm int}(\cK)$.
\end{proof}

\subsection{\bf Relative entropy and matricial distributions }\label{sec:matricialdistributions}

The geometry of convex cones and of the moment problem when $\rho$ is a matricial density function on a compact set $\cS$, as in (\ref{eq:discretemoment}-\ref{generalmomentproblem}), is quite similar to the  case where $\rho$ is only a positive matrix as in Section
\ref{sec:constantdensities}. Appropriate generalizations of the relative entropy functionals allow computable expressions for the corresponding extrema when $\cS$ is a closed interval of the real line, or even a multi-dimensional closed interval in $\mR^k$ ($k>1$). We develop this theory focusing on
(\ref{generalmomentproblem}).

\newcommand{\mmM}{{\mathds{M}}}
\newcommand{\ffM}{{\textswab{M}}}
\newcommand{\ffR}{{\textswab{R}}}
We consider Hermitian $m\times m$ matrix-valued measurable
functions on $\cS$ as a linear space over $\mR$ with an inner product
\[\langle m_1,m_2\rangle =\int_\cS \trace(m_1(\theta)m_2(\theta))d\theta.
\]
We use the notation $\mmM$ to denote 
the Hilbert space of square integrable elements, and the notation $\ffM$ and $\ffM_+$ to denote the cones of elements which are nonnegative and positive definite, respectively, for all $\theta\in\cS$.
The linear operator
\begin{equation}\label{Ltheta}
L\;:\; \ffM \to \fR\;:\; \rho \mapsto R=\int_\cS G_\rleft(\theta)\rho(\theta)G_\rright(\theta)d\theta
\end{equation}
maps $\ffM$ into a subspace of $\mC^{\rleft\times\rright}$ denoted by $\fR$ as before and viewed as a linear space over $\mR$. Both, moments $R$ and their duals $\lambda$ reside in $\fR$ and the geometry is always based on (\ref{innerproduct}). For simplicity of the exposition, we assume that the integration kernels $G_\rleft,G_\rright$ are {\em continuously differentiable} on $\cS$.
The {\em closure} of the range of $\ffM$ is denoted by $\cK=\overline{L(\ffM)}$, while ${\rm int}(\cK)=L(\ffM_+)$. The adjoint transformation is now
\[
L^* : \fR \to \mmM : \lambda \mapsto \rho=\left(G_\rleft(\theta)^* \lambda\, G_\rright(\theta)^*\right)_{\rm Herm}.
\]
It is not difficult to show that the expressions for the dual cone and its interior
\begin{eqnarray*}
\cK^\dual&=&\{\lambda\in\fR\;:\; L^*(\lambda)\in\ffM\}, \mbox{ and }\\
\cK^\dual_+&=&\{\lambda\in\fR\;:\; L^*(\lambda)\in\ffM_+\}
\end{eqnarray*}
remain valid (except for the obvious change where $\ffM$ replaces our earlier $\fM$). The analog of (\ref{nonintersect}) will be needed (in Theorem \ref{thm1:matrix}) which, can also be expressed as
\begin{equation}\label{nonintersect2}
K^\dual_+\neq \emptyset.
\end{equation}
Finally we define as before
\[
\ffM_{R,+}:=\ffM_+\cap \{\rho\in\mmM\;:\; R=L(\rho)\}
\]
as we seek to determine whether or not $\ffM_{R,+}=\emptyset$, or equivalently, whether $R\in{\rm int}(\cK)$.

For future reference we bring in a characterization of elements $R\in\cK$ analogous to the scalar real case given in \cite[page 14]{KreinNudelman}. Given $R\in\fR$, define the real-valued functional
\begin{eqnarray}\nonumber
\fC_R&:&\fR \to \mR\\
&:& \lambda\mapsto \langle\lambda, R\rangle\label{functionalC}
\end{eqnarray}
Such a bounded functional is said to be {\em nonnegative} (resp., {\em positive})---denoted by $\fC_R\geq 0$ (resp., $\fC_R>0$), if and only if the infimum of
$\fC_R(\lambda)$ over $\lambda\in\cK^\dual_+$ of unit norm is positive (resp.\ nonnegative).

\begin{prop}\label{characterizationofR}{\sf The following hold:
\begin{eqnarray*}
R\in\cK&\Leftrightarrow& \fC_R\geq 0\\
R\in{\rm int}(\cK) &\Leftrightarrow&\fC_R>0.
\end{eqnarray*}
}
\end{prop}
\vspace*{.2in}

\begin{proof} We now only prove necessity, which is needed in the proof of Theorem \ref{thm1:matrix}.
The proof of sufficiency will be given at the end Section \ref{minimizersofmbS}.

If $R\in{\rm int}(\cK)$, there exists a particular $\rho\in\ffM_+$ such that $R=L(\rho)$. It readily follows that
$\fC>0$. If $R\in\cK$, there exist 
an approximating sequence $R_i\to R$ ($i=1,2,\dots$)
where $R_i=L(\rho_i)$ and $\rho_i\in\ffM$. If $\lambda\in\cK_+^\dual$ then
\[\fC_R(\lambda)=\lim_{i\to\infty} \langle \lambda,R_i\rangle
\]
which is $>0$ and hence at least $\fC\geq 0$.
\end{proof}

\newcommand{\mbS}{{\mathbbm{S}}}
We now turn to relative entropy functionals for matricial distributions.
Given $\rho,\sigma\in\ffM_+$,
\begin{equation}\label{mbS}
\mbS (\rho\,\|\sigma):=\int_\cS \trace(\rho\,\log\rho-\rho\, \log\sigma) d\theta.
\end{equation}
Once again, minimizers of relative entropy subject to the moment constraints (\ref{generalmomentproblem}) take a particularly simple
form amenable to a numerical solution via continuation methods.
We follow the same plan as in Section \ref{sec:constantdensities} by focusing successively
on each of the two alternative choices, $\mbS(I\|\rho)$ and then $\mbS(\rho\,\|I)$. A significant
departure from the case of constant densities shows up
when considering the dimension of the support set $\cS$ in the context of $\mbS(I\|\rho)$.

\subsubsection{\bf Minimizers of $\mbS(I\|\rho)=-\int_\cS\trace(\log(\rho))d\theta$:}\label{minimizersofmbS}

In complete analogy with constant case the derivative of the Lagrangian
\[
\cL(\lambda,\rho):=-\int_\cS\trace(\log(\rho))d\theta - \langle\lambda,R-L(\rho)\rangle
\]
in the direction $\delta\in\mmM$ is
\begin{eqnarray*}
d \cL(\lambda,\rho \,;\, \delta)&:=&  \trace\int_\cS(-M_{\rho(\theta)}^{-1}+ L^*(\lambda))\delta(\theta) d\theta \\
&=& \trace\int_\cS(-\rho(\theta)^{-1}+ L^*(\lambda))\delta(\theta)d\theta,
\end{eqnarray*}
where, once again, the presence of the $\trace$ allows replacing the ``super-operator'' $M_{\rho(\theta)}^{-1}$
by multiplication by $\rho(\theta)^{-1}$, pointwise over $\cS$.
The fundamental lemma in calculus of variations (see e.g., \cite{bliss})
now gives the stationarity condition
\begin{equation}\label{Lstarinv}
\rho=L^*(\lambda)^{-1}.
\end{equation}
In order for $\rho\in\ffM$ it is necessary that $L^*(\lambda)$ is strictly positive on $\cS$. Thus, we
consider the ``rational'' family of potential minimizers for $\mbS(I\|\rho)$
\begin{eqnarray*}
\ffM_{\rm rat}&:=&\left\{ \rho=L^*(\lambda)^{-1}, \mbox{ with }\lambda\in\cK^\dual_+\right\},
\end{eqnarray*}
where we seek a solution to the moment constraints (\ref{generalmomentproblem}).
It turns out that if a solution exists then, a particular one exists in $\ffM_{\rm rat}$ and that it can be
obtained by computing the fixed point of an exponentially converging matrix differential equation. This differential equation is an appropriate generalization of (\ref{firstoccurence}). We summarize all these conclusions below.

\begin{thm}\label{thm1:matrix}{\sf If ${\rm dim}(\cS)=1$, condition (\ref{nonintersect2}) holds, and $R\in{\rm int}(\cK)$, then $\mbS(I\|\rho)$ has a minimum in $\ffM_{R,+}$ which is unique and belongs to $\ffM_{\rm rat}$. Furthermore, for any  $\lambda_0\in\cK_+^\dual$, the solution $\lambda_t$ of
the matrix differential equation
\begin{eqnarray}\label{diffeqlambda}
\frac{d}{dt} \lambda_t = \left( \nabla h|_{\lambda_t}\right)^{-1} (R-L(L^*(\lambda_t)^{-1})),
\end{eqnarray}
where
\begin{equation}\label{nablathm}
\nabla h|_{\lambda_t}:\fR\to\fR\;:\; \delta \mapsto 
L(L^*(\lambda_t)^{-1}L^*(\delta)L^*(\lambda_t)^{-1}),
\end{equation}
belongs to $\cK_+^\dual$ for all $t\in[0,\infty)$, it converges
to a point $\hat{\lambda}\in\cK_+^\dual$ as $t\to\infty$ corresponding to this unique minimizer
$\rho=L^*(\hat{\lambda})^{-1}$ for $\mbS(I\|\rho)$ satisfying
$R=L(\rho)$. The differential equation (\ref{diffeqlambda}) is exponentially convergent as the square distance
$V(\lambda_t)=\|R-L(L^*(\lambda_t)^{-1})\|^2$
satisfies
\[
\frac{dV(\lambda_t)}{dt}=-2V(\lambda_t).
\]
Conversely, if $R\not\in{\rm int}(\cK)$ and the dimension of $\cS$ is one, then the differential equation (\ref{diffeqlambda}) diverges.}
\end{thm}

Equations (\ref{diffeqlambda}) is equivalent to
\begin{eqnarray}\label{diffeqlambda_alt}
\frac{d}{d\tau} \lambda_\tau= \left( \nabla h|_{\lambda_\tau}\right)^{-1} (R-R_0),
\end{eqnarray}
modulo scaling of the integration variable (see below).
The latter can be integrated over $[0,1]$, and then $\hat{\lambda}=\lambda_\tau|_{\tau=1}$, yet (\ref{diffeqlambda}) appears
preferable for numerical reasons.

We wish to point out that the assumption on the dimension of $\cS$ can be slightly relaxed to being at most two provided $\cS$ is a torus and $G_\rleft,G_\rright$ doubly periodic accordingly, (cf.\ \cite[Example 2 on page 882]{LangMcClellan}, \cite{ac_may2004}).

\begin{proof}
Once again we consider the mapping
\[
h\;:\;\fR\to\fR \;:\; \lambda\mapsto R=L(L^*(\lambda)^{-1}).
\]
The Jacobian $\nabla h|_\lambda$ is Hermitian and invertible when $\lambda\in\cK^\dual_+$. The proof is identical to the one given in Section \ref{sec:minimizersofrho}.
Choose $\lambda_0\in\cK^\dual_+$, set $R_0=L(L^*(\lambda_0)^{-1})$, and
consider the one-parameter homotopy of maps
\begin{equation}\label{homotopy}
L(L^*(\lambda_\tau)^{-1}) = R+\tau(R-R_0) =:R_\tau,
\end{equation}
for $\tau\in[0,1]$. The idea is to follow a path of solutions $\lambda_\tau$ ($\tau\in[0,1]$), ensure that
this is contained in $\cK_+^\dual$, and set $\hat{\lambda}=\lambda_1$ which then satisfies the sought conditions ($L^*(\hat{\lambda})^{-1}\in\ffM_{\rm rat}$ and $L(L^*(\hat{\lambda})^{-1})=R$).

Clearly, $L^*(\lambda_0)^{-1}\in\ffM_{\rm rat}$ and $R_0\in{\rm int}(\cK)$. If $R\in {\rm int}(\cK)$ then so is $R_\tau$ for $\tau\in[0,1]$. We claim that for all $\tau\in[0,1]$ there exists $\lambda_\tau\in\cK_+^\dual$ such that
\begin{equation}\label{Rlambdatau}R_\tau=L(L^*(\lambda_\tau)^{-1})
\end{equation}
as before. 
The arguments are similar to those given in the proof of Theorem \ref{thm1a} and are based on
the fact that $h$ is a local diffeomorphism. Values $\lambda_\tau$ which obey
\begin{equation}\label{unscaled}
\frac{d\lambda_\tau}{d\tau}=(\nabla h_{\lambda_\tau})^{-1}(R-R_0),
\end{equation}
satisfy (\ref{Rlambdatau}) as long as the path stays in $\cK^\dual_+$. We need to rule out  $\lambda_\tau$ crossing the boundary of $\cK^\dual_+$ or tending to $\infty$ at a $\tau\leq 1$. Either possibility contradicts $R_\tau\in{\rm int}(K)$ ($\tau\in[0,1]$) in a way analogous to the earlier arguments in the proof of Theorem \ref{thm1a}.

We consider a maximal interval $[0,\epsilon)$ over which $R_\tau\in{\rm int}(K)$,
and note that
 \[
 \langle\lambda_\tau,R_\tau\rangle=\int_\cS\trace(I)d\theta= \trace(I)\cdot{\rm measure}(\cS) <\infty,
 \]
on $[o,\epsilon)$. If  $\|\lambda_\tau\|\to\infty$, then $\fC_{R_\tau}$ cannot be bounded away from zero 
since $ \langle\frac{\lambda_\tau}{\|\lambda_\tau\|},R_\tau\rangle\to 0$. Hence, $R_\epsilon\not\in{\rm in}(\cK)$ and $\epsilon>1$.

If $\lambda_\epsilon$ lies on the boundary of $\cK^\dual$, then $L^*(\lambda_\epsilon)$ has a root in $\cS$. It is here that the dimension of $\cS$ becomes important. As long as the dimension of $\cS$ is one, $L^*(\lambda_\epsilon)^{-1}$ is not integrable as a function of $\theta\in\cS$ and
$R_\epsilon=L(L^*(\lambda_\epsilon)^{-1})$ cannot be finite. Thus again, $\epsilon>1$.

We finally express (\ref{unscaled}) in a ``feedback form''.
We first replace $\tau$ with $t=-\log(1-\tau)$. In this case, $\tau=1-e^{-t}$ and $t$ varies in $[0,\infty)$ as $\tau$ varies in $[0,1]$. If we denote $\lambda_t:= \lambda_{\tau(t)}$ and $R_t:=R_{\tau(t)}$, then
\begin{eqnarray}\nonumber
\frac{d}{dt} \lambda_t &=&\left(\frac{d\tau}{dt}\right)\left( \nabla h|_{\lambda_t}\right)^{-1} (R-R_0)\\
&=& (1-\tau (t))\left( \nabla h|_{\lambda_t}\right)^{-1} (R-R_0)\nonumber\\
&=& \left( \nabla h|_{\lambda_t}\right)^{-1} (R-R_t).\label{diffeqfinal}
\end{eqnarray}
The same substitution gives
$\frac{dR_t}{dt}=R_1-R_t$,
and that
\[
V(\lambda_t):=\|R-R_t\|^2=\|R-L(L^*(\lambda_t)^{-1})\|^2
\]
satisfies
\begin{eqnarray*}
\frac{dV(\lambda_t )}{dt}&=&-2\langle R-L(L^*(\lambda_t)^{-1}),\nabla h|_{\lambda_t} \frac{d}{dt}\lambda_t\rangle\\
&=&-2\langle R-R_t,R-R_t\rangle=-2V(\lambda_t).
\end{eqnarray*}

Uniqueness of the representation $R=L(\rho)$ with $\rho\in\ffM_{\rm rat}$
follows from the fact that such a $\rho$ is a minimizer of the convex functional $\mbS(I\|\rho)$ subject to the moment constraints. An alternative but equivalent argument can be based on the fact that
$h$ is a $C^1$-mapping between open convex subsets of a Euclidean space, with a positive definite Jacobian everywhere.

In the other direction, if $R\not\in {\rm int}(\cK)$, then the path $R_\tau$
either crosses or at least, in case $R\in\partial\cK$, tends to the boundary of  $\cK$.
In either case, $\lambda_t$ grows unbounded and the differential equation diverges.
\end{proof}

We now complete the proof of Proposition \ref{characterizationofR}.

\begin{proof}{\em $[$``sufficiency'' in Proposition \ref{characterizationofR}$\;]\;\;$}
In the proof of Theorem \ref{thm1:matrix}, in essence, we used the positivity of $\fC_{R}$ to obtain a representation $R=L(\rho)$ with $\rho\in\ffM_+$. Thus, the same line of argument gives that if
$\fC_{R}>0$ then $R\in{\rm int}(\cK)$.

If $\fC_R\geq 0$ but not necessarily $>0$, then chose $R_0\in{\rm int}(\cK)$ and note that $\fC_{R+\frac{1}{k}R_0}>0$ for $k=1,2,\ldots$. Obviously, $R+\frac{1}{k} R_0\in{\rm int}(\cK)$ ($k\to \infty$) and $R$ is at least in $\cK$.\footnote{A uniform bound on
the integral of densities corresponding to $R+\frac{1}{k} R_0$ can be shown.
This can be used to establish a finite nonnegative measure $d\mu$ corresponding to $R$.}
\end{proof}

\subsubsection{\bf Minimizers of $\mbS(\rho\|I)=\int_\cS\trace(\rho\log(\rho))d\theta$:}\label{exp}

Once again, the derivative of the Lagrangian
\[
\cL(\lambda,\rho):=-\int_\cS\trace(\rho\log(\rho))d\theta - \langle\lambda,R-L(\rho)\rangle
\]
in the direction $\delta\in\mmM$ is
\begin{eqnarray*}
d \cL(\lambda,\rho \,;\, \delta)&:=&  \trace\int_\cS(\log(\rho)+I+ L^*(\lambda))\delta(\theta) d\theta.
\end{eqnarray*}
The stationarity condition leads to the expression
\begin{equation}\label{expminusLstar}
\rho=\frac{1}{e}\exp(-L^*(\lambda)),
\end{equation}
for the minimizer, except that now $\rho$ is a function of $\theta\in\cS$. We
consider the ``exponential'' family
\begin{eqnarray*}
\ffM_{\rm exp}&:=&\left\{ \rho=\frac{1}{e}\exp(-L^*(\lambda)), \mbox{ with }\lambda\in\fR\right\},
\end{eqnarray*}
of potential minimizers for $\mbS(\rho\,\|I)$, where we seek a solution to (\ref{generalmomentproblem}). The development runs in parallel to the case where $\rho\in\ffM_{\rm rat}$ with one important difference. The ``Lagrange multipliers'' $\lambda$ no longer need to be restricted to $\cK^\dual$ and
existence of solutions when $R\in\cK$ can be guaranteed even when ${\rm dim}(\cS)>1$. Moreover,
(\ref{nonintersect2}) is no longer necessary and existence of solution to the moment problem in $\ffM_{\rm exp}$ is impervious to the dimension of the dual cone $\cK^\dual_+$.

\begin{thm}\label{thm2:matrix}{\sf If $R\in{\rm int}(\cK)$ then the entropy functional $\mbS(I\|\rho)$ has a minimum in $\ffM_{R,+}$ which is unique and belongs to $\ffM_{\rm exp}$. Furthermore, for any  $\lambda_0\in\fR$, the solution $\lambda_t$ of
\begin{eqnarray}\label{diffeqlambda2}
\frac{d}{dt} \lambda_t = \left( \nabla k|_{\lambda_t}\right)^{-1} (R-L(\lambda_t)),
\end{eqnarray}
where
\begin{equation}\label{nablathm2}
\nabla k|_{\lambda_t}:\fR\to\fR\;:\; \delta \mapsto 
-\frac{1}{e}L(M_{\exp(-L^*(\lambda_t))}(L^*(\delta)),
\end{equation}
remains bounded for $t\in[0,\infty)$ and converges
to $\hat{\lambda}\in\fR$ as $t\to\infty$ corresponding to the unique minimizer
$\rho=\frac{1}{e}\exp(-L^*(\hat{\lambda}))$
for $\mbS(I\|\rho)$ subject
to $R=L(\rho)$.
The convergence is exponential as
$V(\lambda_t)=\|R-\frac{1}{e}L(\exp(-L^*(\hat{\lambda_t})))\|^2$
satisfies
$\frac{dV(\lambda_t)}{dt}=-2V(\lambda_t)$.
Conversely, if $R\not\in{\rm int}(\cK)$ then the differential equation (\ref{diffeqlambda2}) diverges.}
\end{thm}

\begin{proof} 
The arguments are for the most part identical to those used in proving Theorem \ref{thm2a},
i.e., we now consider
\[
\kappa\;:\;\fR\to{\rm int}(\cK)\subset \fR \;:\; \lambda\mapsto R=L(\frac{1}{e}e^{-L^*(\lambda)}),
\]
observe that the Jacobian is negative definite for any value of $\lambda$, and use it to track
a linear path $(1-\tau)R_0+\tau R$ ($\tau\in[0,1]$) from $R_0=L(\frac{1}{e}e^{-L^*(\lambda_0)})$ to the given $R$ in the $\lambda$-coordinates. This is done by integrating (\ref{firstoccurenceexp}) over $[0,1]$
starting from arbitrarily chosen starting point $\lambda_0$. By construction, the solution of (\ref{firstoccurenceexp}) corresponds
to $\rho_\tau=\frac{1}{e}e^{-L^*(\lambda_\tau)}\in\ffM_+$ which satisfies $R_\tau=L(\rho_\tau)$.
It is clear that if $R\not\in{\rm int}(\cK)$, then the differential equation diverges for $\tau\leq 1$ (since otherwise $\rho_\tau\in\ffM_+$ and $R_\tau=L(\rho_\tau)$ would hold on $[0,1]$, contradicting $R\not\in{\rm int}(\cK)$). 
We only need to show that if $R\in{\rm int}(\cK)$, then $\lambda_\tau$ remains bounded
for $\tau\in[0,1]$. This yields $\hat{\lambda}=\lambda_\tau|_{\tau=1}$ which satisfies $R=L(\frac{1}{e}e^{-L^*(\hat{\lambda})})$. Then (\ref{diffeqlambda2}) can be obtained via a change of variables as in
Theorem \ref{thm1:matrix}. The same applies to deriving the differential equation for the ``error'' $V(\lambda_t)$.

In order to show that, in the event $R\in{\rm int}(\cK)$, $\lambda_\tau$ remains bounded on $[0,1]$ we extend the argument used to prove Theorem \ref{thm2a} to the present case where $\rho$ is a matrix valued function on $\cS$. The key is to observe that, when (\ref{firstoccurenceexp}) is integrated over a maximal interval $[0,\epsilon)$, any convergent subsequence of $\ell_\tau:=\lambda_\tau/\|\lambda_\tau\|$ ($\tau\in[0,\epsilon)$) must have a limit point $\ell_\epsilon$ for which $L^*(\ell_\epsilon)\in\ffM$ but not in $\ffM_+$. Moreover, $L^*(\ell_\epsilon)$ must be singular on $\cS_0\subset \cS$ (a subset of possibly zero measure).
To see this note the following.
If $L^*(\ell_\epsilon)\in\ffM_+$ then $L^*(\ell_{\tau_i})$  is bounded away from zero and positive for $i$ large enough, whereas if $L^*(\ell_\epsilon)\not\in\ffM$ then there is a subset of $\cS$
of non-zero measure where $L^*(\ell_{\tau_i})$ is negative. Either way
$L(\exp(-L^*(\lambda_{\tau_i})))= L(\exp(-L^*(\ell_{\tau_i}\cdot\|\lambda_{\tau_i}\|)))$ cannot tend to 
$R_\epsilon$ as it should. In the first instance it goes to zero and in the second it becomes unbounded.
Thus $L^*(\ell_\epsilon)\in\ffM$ but singular for certain values of $\theta$.
Below we show that this implies $R_\epsilon\not\in{\rm int}(\cK)$, which then proves that $\epsilon>1$ and that (\ref{diffeqlambda2}) can be integrated on $[0,1]$.

To show that $R_\epsilon\not\in{\rm int}(\cK)$ it suffices to show that $\fC_{R_\epsilon}$ is not strictly positive. To this end, we evaluate
\begin{eqnarray*}
\fC_{R_{\tau_i}}(\ell_\epsilon)&=&
\langle \ell_\epsilon,L(\exp(-L^*(\ell_{\tau_i}\|\lambda_{\tau_i}\|)))\rangle\\
&=& \langle L^*(\ell_\epsilon),\exp(-L^*(\ell_{\tau_i}\|\lambda_{\tau_i}\|)))\rangle\\
&=&\int_\cS \trace (L^*(\ell_\epsilon)\exp(-L^*(\ell_{\tau_i}\cdot \|\lambda_{\tau_i}\|))) d\theta\\
&=&\int_\cS \trace (L^*(\ell_\epsilon)(\exp(-L^*(\ell_{\tau_i})))^{\|\lambda_{\tau_i}\|} d\theta.
\end{eqnarray*}
For each value of $\theta$, $(\exp(-L^*(\ell_{\tau_i})))^{\|\lambda_{\tau_i}\|}$ tends to zero outside
the null space of $L^*(\ell_\epsilon)$. Since $\cS$ is compact,
the integrand goes to zero uniformly in $\theta$ as $i\to\infty$. Therefore, $\fC_{R_{\tau_i}}(\ell_\epsilon)\to 0$ as well, and $R\not\in{\rm int}(\cK)$.\end{proof}

\subsection{\bf Non-equispaced arrays (cont.)}

\begin{figure}[htb]\begin{center}
\includegraphics[totalheight=7cm]{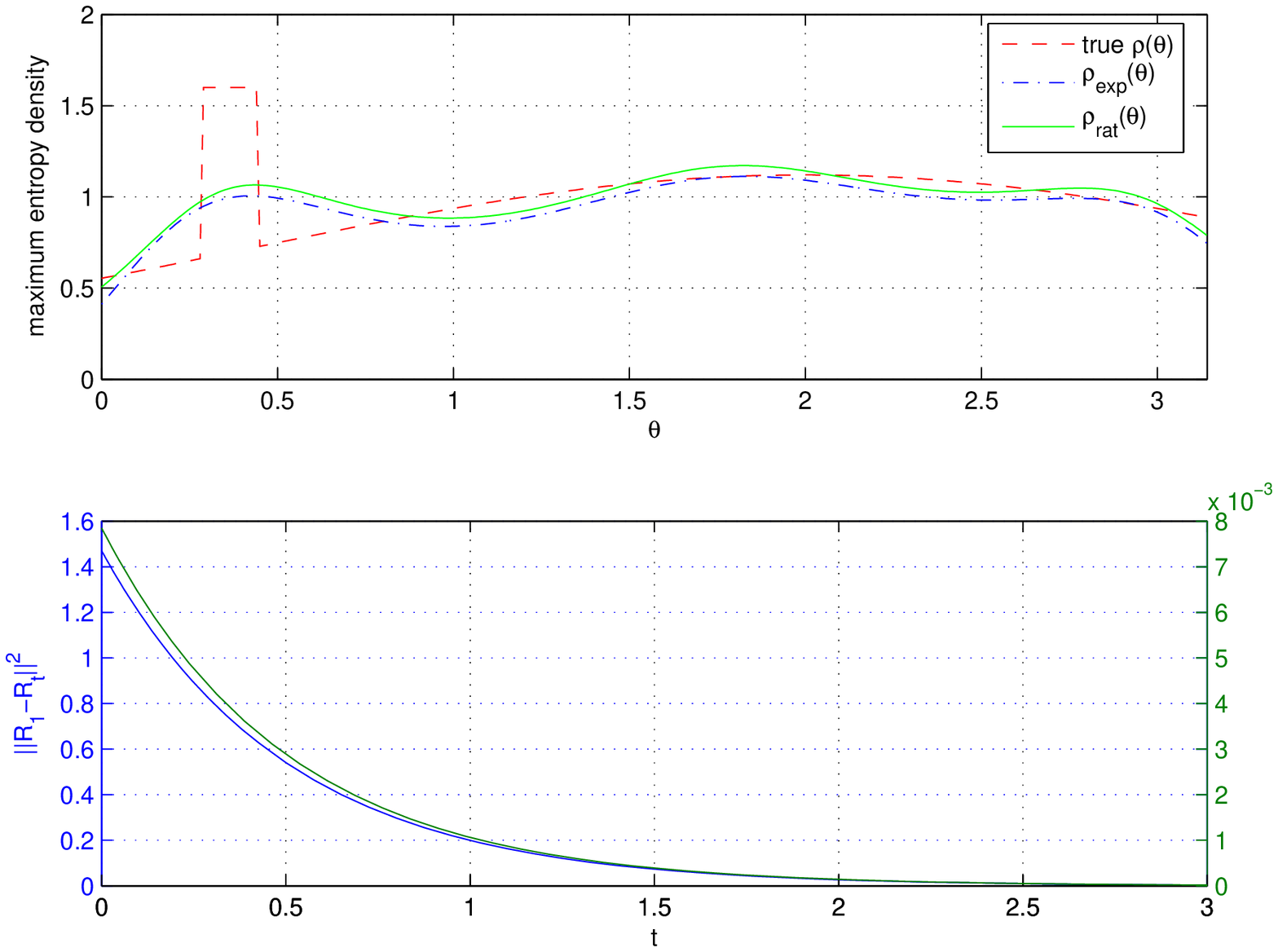}
\caption{}\label{fig:plot1}\end{center} \end{figure}

We continue with Example \ref{example1}.
We begin with  a ``true'' density $\rho_{\rm true}$ shown in Figure \ref{fig:plot1}
and generate covariance samples $R$. This ``true''
density does not need to be in any particular form---computation of $R$ is done via numerical integration.

Next, we integrate (\ref{diffeqlambda}) and (\ref{diffeqlambda2})
taking $\lambda_0=\left[\begin{matrix}1&0&0&0\end{matrix}\right]$, and display in
Figure \ref{fig:plot1} the resulting $\rho_{\rm  exp}(\lambda_\infty,\theta)$ and $\rho_{\rm  rat}(\lambda_\infty,\theta)$, for comparison. Both are constructed using the fixed point of the corresponding differential equations. The rate of convergence is the same, while the distance of the starting choice (for the same $\lambda_0$) may be different---as is the case here (with $\ffM_{\rm exp}$ corresponding to the $y$-axis to the left and $\ffM_{\rm rat}$ the labeling to the right in subplot $(2,1)$).

\section{\bf The complete set of positive solutions}\label{completesolutions}

Reference \cite{GL} suggested that all positive solutions to the moment 
problem may be obtained as minimizers of a suitable entropy functional, e.g., as being
\begin{equation}\label{weightedentropy}
{\rm argmin}\{\mS(\sigma\,\|\rho)\;:\; R=L(\rho)\}
\end{equation}
with $\sigma$ thought of as a parameter. This was carried out successfully in \cite{GL} and \cite{ac_may2004} for the case where density functions are scalar-valued, for different levels of generality. Naturally, certain complications arise in the matricial setting. We discuss this next in the context of constant $\rho,\sigma$ as in Section \ref{sec:constantdensities}). The generalization to the non-constant case is straightforward and a positive result is given for the general case.

Considering the Lagrangian and the stationarity conditions for (\ref{weightedentropy}) we arive at
\begin{eqnarray*}
d \cL(\lambda,\rho\,;\,\delta)&=&
 \trace(-\delta M_\rho^{-1} (\sigma))+ \langle L^*(\lambda),\delta\rangle\\
&=& \trace\left(-\delta M_\rho^{-1} (\sigma) + \delta L^*(\lambda)\right),
\end{eqnarray*}
leading to
\[
M_\rho^{-1}(\sigma) = L^*(\lambda).
\]
Although the ``parameter'' $\sigma$ can be readily
expressed as $M_\rho(L^*(\lambda))$, the density $\rho$ which we are interested in, cannot be
expressed in any effective way as a function of $\sigma$ and the dual variable $\lambda$. Thus, a convenient functional form for the minimizer of (\ref{weightedentropy}) is unkown.

The option of minimizing
$\mS(\rho\,\|\sigma)$ subject to $R=L(\rho)$
however, goes through. Analysis of the corresponding Lagrangian readily leads to
\[
\rho=\frac{1}{e}\exp(\log(\sigma)-L^*(\lambda)).
\]
A computational theory, following the lines of Sections \ref{expconstant} and \ref{exp} easily carries through.

An attractive third alternative originates in the observation that the geometry of the problem,
throughout, was inherited by the definiteness of the Jacobian maps.
This suggests to forgo an explicit form for the entropy functional and start instead with a computable
Jacobian. To this end we consider
\begin{eqnarray*}
h_\sigma&:&\lambda\mapsto L(\sigma^{1/2} L^*(\lambda)^{-1} \sigma^{1/2} ), \mbox{ and}\\
\kappa_\sigma&:&\lambda\mapsto L(\sigma^{1/2} \frac{1}{e}\exp(-L^*(\lambda))\sigma^{1/2}).
\end{eqnarray*}
The respective Jacobians are
\begin{eqnarray*}
\nabla h_\sigma|_\lambda&:&\delta \mapsto L(\sigma^{1/2} L^*(\lambda)^{-1} L^*(\delta)L^*(\lambda)^{-1}\sigma^{1/2} ), \mbox{ and}\\
\nabla \kappa_\sigma|_\lambda&:&\delta \mapsto \frac{1}{e}L(\sigma^{1/2}M_{\exp(-L^*(\lambda))}(-L^*(\delta))\sigma^{1/2}).
\end{eqnarray*}
They are both sign definite as before and, almost verbatim, we can replicate the conclusions of Theorems \ref{thm1:matrix} and \ref{thm2:matrix}. These are combined into the following statement.

\begin{thm}\label{thm1:matrixW}{\sf Let $R\in{\rm int}(\cK)$ and $\sigma\in\ffM_+$.
If ${\rm dim}(\cS)=1$, condition (\ref{nonintersect2}) holds, and $\lambda_0\in{\rm int}(\cK_+^*)$, then the solution to
\begin{eqnarray}\label{diffeqlambdaW}
\frac{d}{dt} \lambda_t = \left( \nabla h_\sigma |_{\lambda_t}\right)^{-1} (R-h_\sigma(\lambda_t))
\end{eqnarray}
remains in $\cK_+^\dual$ for $t\geq 0$ and as $t\to\infty$ converges to a unique value $\lambda_{\rm r}\in\cK_+^\dual$ such that
$R=h_\sigma(\lambda_{\rm r})$.
On the other hand, for any $\lambda_0\in\fR$ the solution to
\begin{eqnarray}\label{diffeqlambdaW2}
\frac{d}{dt} \lambda_t = \left( \nabla \kappa_\sigma |_{\lambda_t}\right)^{-1} (R-\kappa_\sigma(\lambda_t))
\end{eqnarray}
remains bounded for $t\geq 0$ and as $t\to\infty$ converges to a unique value ${\lambda_{\rm e}}\in\cK_+^\dual$ such that
$R=\kappa_\sigma(\lambda_{\rm e})$.
In case $R\not\in{\rm int}(\cK)$, then (\ref{diffeqlambdaW2}) diverges.
In case $R\not\in{\rm int}(\cK)$ and ${\rm dim}(\cS)=1$, then (\ref{diffeqlambdaW}) diverges as well. 
}
\end{thm}

The importance of recasting Theorems \ref{thm1:matrix} and \ref{thm2:matrix}
as above, by incorporating arbitrary $\sigma$'s in $\ffM_+$, allows obtaining {\em any} density function which is consistent with the data $R$ by such a procedure. To see this note that, if $\rho$ consistent with the data, then working backwards we can select $\sigma$ accordingly so that
$\rho$ equals $\sigma^{1/2} L^*(\lambda)^{-1} \sigma^{1/2} $ or $\frac{1}{e}\sigma^{1/2}\exp(-L^*(\lambda))\sigma^{1/2}$ for any $\lambda$ (in $\cK^\dual_+$ and $\fR$, respectively).
Thus, Theorem \ref{thm1:matrixW} gives descriptions of {\em all} positive densities that are consistent with the data $R$---simply choose the ``correct'' $\sigma$.
  
A potentially important application is when prior information may dictate a choice of $\sigma$.
In this case, using Theorem \ref{thm1:matrixW} we may obtain an admissible
density function which is ``closer to our expectations.'' We amplify this remark by reworking Example \ref{example1} with a suitable weight.

\subsection{\bf Non-equispaced arrays (cont.)}\label{incorporate}

\begin{figure}[htb]\begin{center}
\includegraphics[totalheight=7cm]{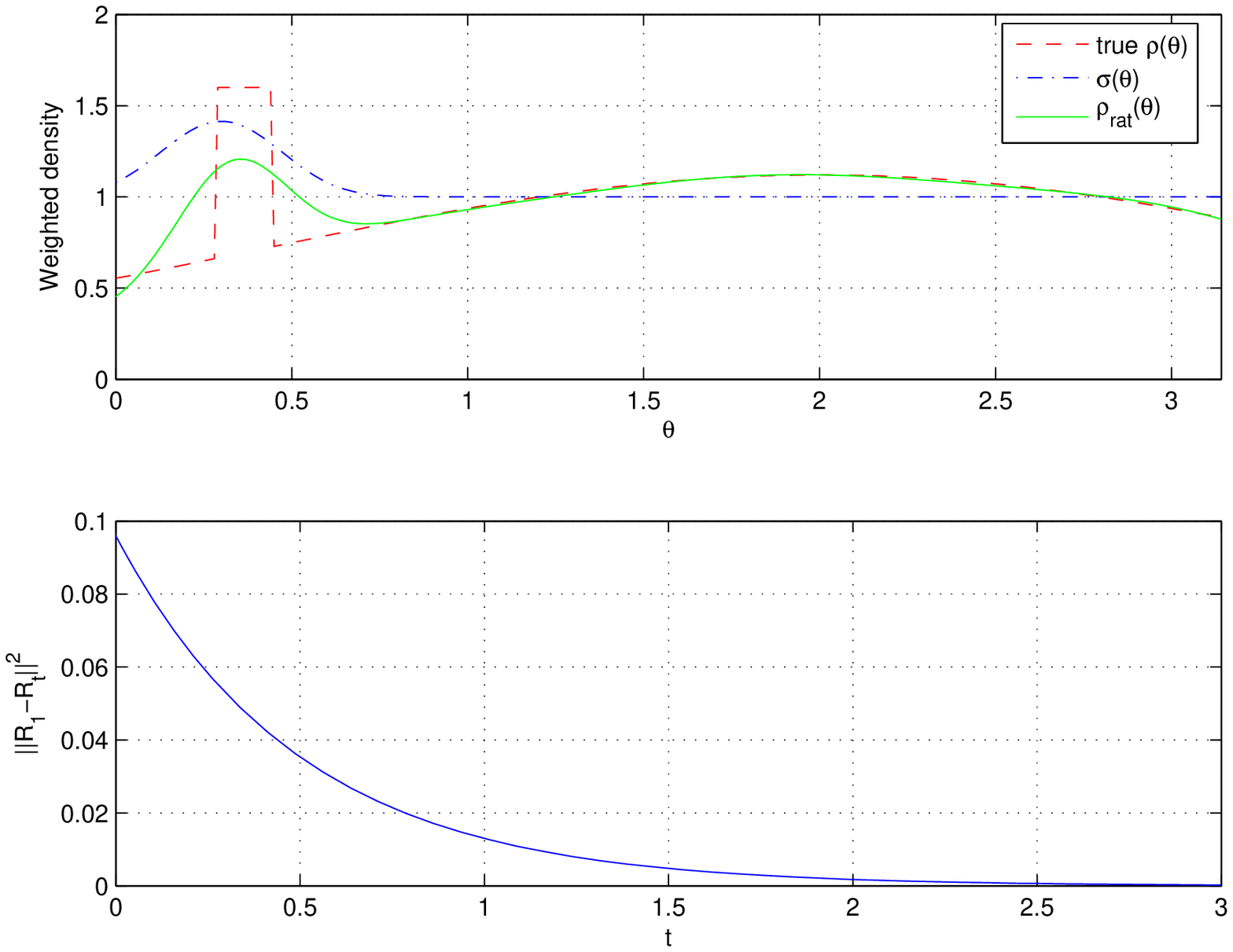}
\caption{}\label{fig:plot1W}\end{center} \end{figure}

Figure \ref{fig:plot1W} compares the ``true'' density function $\rho_{\rm true}$ which was used to generate the moments,
and a density  $\rho_{\rm rat}=\sigma^{1/2}L^*(\lambda)\sigma^{1/2}$ which is computed according to Theorem \ref{thm1:matrixW}. The original density $\rho_{\rm true}$ has discontinuous peak at about $\theta\sim 0.35$. Then $\sigma$ has been selected so as to be $\geq 1$ in the neighborhood of $\theta\sim 0.35$---actually centered about $0.25$. (The accuracy of the ``match'' does not seem critical.)
The density $\rho_{\rm rat}$ is seen to be a better match as compared with the ``unweighted'' case of Figure \ref{fig:plot1}.
Subplot (2,1) shows the value of $\|R_1-R_t\|$ as before, and highlights the fact that, again,
$\rho_{\rm rat}$ is consistent with the moments. Since $\lambda_0 G(\theta)\equiv 1$, if we choose $\sigma=\rho_{\rm true}$ (using $100\%$ hindsight), we obviously obtain a perfect match as explained above.

\subsection{\bf State-statistics and analytic interpolation with degree constraint (cont.)}\label{example4_continued}

In Section \ref{completesolutions}, the map $h_\sigma$ can be replaced by
\[
\lambda\mapsto L(\varphi L^*(\lambda)^{-1} \varphi^* )
\]
where $\sigma=\varphi \varphi^*$ is a factorization of $\sigma$ with $\varphi$ not necessarily Hermitian, with the obvious modifications in the expression for the corresponding Jacobian. The statement of the theorem holds with no changes. The same applies to $\kappa_\sigma$ which can also be cast with respect to an arbitrary factorization of $\sigma$---but this will not concern us here. Instead, we consider the  setting of Section \ref{example4} where
\[
L^* \;:\; \lambda \mapsto B^*(I-e^{-j\theta} A^*)^{-1}\lambda(I-e^{j\theta}A)^{-1}B.
\]
If we take $\varphi(z)=I+C_o z (I-zA)^{-1}B$ so that $\varphi^{-1}$ is also analytic in $\mathbb D$ (which corresponds to $C_o$ chosen so that $A-BC_o$ is a Hurwitz matrix), then the resulting density function
\begin{eqnarray*}
\rho(\theta)&=&\varphi L^*(\lambda)^{-1} \varphi^*\\
&=& (G_o(e^{j\theta})^* \lambda G_o(e^{j\theta}))^{-1},
\end{eqnarray*}
with $G_o(z)=(I-z(A-BC_o))^{-1}B$.
This is a rational spectral density of degree at most twice the dimension of (\ref{i2s}), and hence, it gives rise to a positive-real interpolant $F$ as in (\ref{correspondingpr}) of McMillan degree at most equal to the dimension of (\ref{i2s}). 

\section{\bf Concluding remarks}\label{conclusion}

We presented an approach for constructing matrix-valued density functions which are consistent with given moments. Section \ref{completesolutions} describes, in the spirit of the mathematical theory on the moment problem,  all positive-definite density functions which are consistent with the data. The non-parametric description given in Section \ref{completesolutions} (non-parametric since it amounts to an arbitrary choice of a weight-density $\sigma$) should prove useful in case we wish to incorporate prior information (e.g., subsection \ref{incorporate}, and cf.\ \cite{GL}).

The basic problem of characterizing admissibilitiy of  a matricial moment $R$ has been cast in terms of the positivity of a suitable functional, $\fC_R$, in complete analogy with the classical case \cite{KreinNudelman}. However, testing for positivity of such a functional is not a trivial matter. In the classical theory, the ``shift'' structure of the space of integration kernels (\cite{KreinNudelman}, \cite{ShohatTamarkin}, \cite{Akhiezer}, \cite{AkhiezerKrein}) allows a simple description of all positive elements in their span, via ``sums of squares.'' This is not the case here. Instead, we determine admissibility of $R$ from the convergence of the differential equation given in e.g., Theorems \ref{thm1:matrix}, \ref{thm2:matrix}.
Yet, a more direct analog of the Pick operator and a corresponding test that would allow a ``certificate of positivity of $\fC_R$,'' would be highly desirable.

The present work has been influenced by recent literature on ``moments with complexity constraints'' \cite{BGL1,BGL2,Blomqvist,Enqvist,ac,siam,megretski}, and in particular by Byrnes, Gusev and Lindquist \cite{BGuL} who first exploited entropy functionals in such a context. Interpolation or moment problems with degree constraint seek to parametrize solutions of bounded degree within the ``rational familiy.'' The ``trigonometric moment problem with degree constraint'' was first studied in \cite{thesis} in both the scalar and the multivariable setting (via degree theory and homotopy \cite[page 76, and Chapters IV and V]{thesis}). All subsequent literature on ``complexity constraints'' focused on scalar problems until Blomqvist etal.\ \cite{Blomqvist} who study matricial Nevanlinna-Pick interpolation via minimizers of an entropy functional. The framework of the present work, which is also based on entropy functionals, when specialized to analytic interpolation, allows dealing with the most general tangential (and bi-tangential, cf.\ \cite{Arov,DD}) Carath\'{e}odory-Nevanlinna-Pick problems with degree constraint. ``Tangential interpolation'' refers to the case $V(z)$ in e.g., (\ref{PRNehari}-\ref{V}), is a matricial inner factor as opposed to simply a scalar-inner times the identity.  This is examplified in Section \ref{example4_continued}.
While the main focus of the present work remains the general moment problem, consequences of the theory as in Section  \ref{example4_continued} should prove useful in multivariable feedback design with degree constraints \cite{GL2}.

\section*{Acknowledgment}
The author would like to thank Pablo Parillo for his input, and Laurent Baratchart for helpful discussions during a 2003 stay at the Mittag-Leffler Institute in Sweden.

\section{\bf Appendix: Matrix calculus}\label{appendix}

We assemble a number of basic mathematical formulae.
These are expressions for the differential of the matrix exponential and the matrix logarithm
that have been used in the physics literature and in quantum information theory.

\subsection{\bf The matrix exponential}
We begin with the differential of the matrix exponential (see \cite{feynman,najfeld}).
Following \cite[page 164]{heims}, integrate both sides of
\[
\frac{d}{dt}\left[e^{-tA}e^{t(A+\epsilon B}\right]=e^{-tA}\epsilon B e^{t(A+\epsilon B)}
\]
between $0$ and $t$ to obtain
\begin{eqnarray*}
e^{-tA}e^{t(A+\epsilon B)} - I&=&\int_0^te^{-t_1A}\epsilon B e^{t_1(A+\epsilon B)} dt_1.
\end{eqnarray*}
Then
\begin{eqnarray*}
e^{t(A+\epsilon B)} &=& e^{tA}+e^{tA} \int_0^te^{-t_1 A}\epsilon B e^{t_1 (A+\epsilon B)} dt_1\\
&=& e^{tA}+e^{tA} \int_0^te^{-t_1 A}\epsilon B \times \\
&& \hspace*{-25pt} \times \left(e^{t_1A}+e^{t_1A} \int_0^{t_1}e^{-t_2 A}\epsilon B e^{t_2 (A+\epsilon B)} dt_2\right)
dt_1\\
&=& S_0(t) + \epsilon S_1(t) + \epsilon^2 S_2(t) +\ldots
\end{eqnarray*}
where
\begin{eqnarray}
S_0(t) &=&e^{tA} \label{series}\\
S_1(t) &=&e^{tA} \int_0^t e^{-t_1 A} B e^{t_1 A} dt_1  \nonumber \\
S_2(t) &=&e^{tA} \int_0^t \int_0^{t_1} e^{-t_1 A} B e^{(t_1-t_2) A} B e^{t_2 A} dt_2 dt_1, \nonumber 
\end{eqnarray}
and the general term $S_n(t)$ is
\[
e^{tA} \int_0^t\int_0^{t_1}\ldots \int_0^{t_{n-1}}
(e^{-t_1 A} B e^{t_1 A})
\ldots
(e^{-t_n A}B e^{t_n A})
dt_n \ldots dt_1.
\]
We are only interested in the first two terms of this convergent series.

Evaluating at $t=1$ and replacing $\epsilon B$ by $\Delta$, we obtain
\[
e^{A+\Delta} - e^{A} = \int_0^1e^{(1-\tau) A} \Delta e^{\tau A} d\tau + o(\|\Delta\|).
\]
Hence, the differential in the direction $\Delta$
(often refered to as Gateaux, or polar, or Fr\'{e}chet) is given by the linear map
\[
\Delta \mapsto \int_0^1e^{(1-\tau) A} \Delta e^{\tau A} d\tau.
\]
This map represents a ``scrambled'' multiplication of $\Delta$ by $e^{A}$. To see this assume
that $A$ and $\Delta$ commute. Then the right hand side becomes simple $e^{A}\Delta$.

The $S_2$-term in (\ref{series}) gives
the quadratic term in $\Delta$ in the expansion of $e^{A+\Delta} - e^{A}$ as
\[
\Delta \mapsto \int_0^1\left( e^{(1-\tau_1)A}\Delta e^{\tau_1 A} \int_0^1 \left( e^{-A\tau_1 \tau_2}\Delta e^{\tau_1\tau_2A}\right) d\tau_2  \right) \tau_1d\tau_1.
\]

In general, for Hermitian matrices $C$ and $\Delta$, and $C\geq 0$, define the ``non-commutative'' or ``scrambled'' multiplication of $\Delta$ by $C$ via the operator
\begin{equation}\label{multiplicationoperator}
M_C\;:\; \Delta \mapsto \int_0^1C^{(1-\tau)} \Delta C^{\tau} d\tau.
\end{equation}
This gives a compact expression for the differential of $e^A$, summarized below.

\begin{prop}
{\sf The differential of $\exp(A):=e^A$ is $M_{e^A}$.}
\end{prop}

\subsection{\bf The matrix logarithm}
We now turn to the matrix logarithm. Integrate both sides of
\[
\frac{d}{dt}\left[ \log(I+tP)-\log(I+tQ)\right]=(I+tP)^{-1}P-Q(I+tQ)^{-1}
\]
between $0$ and $1$ to obtain that
\[
\log(I+P)-\log(I+Q)=\int_0^1(I+tQ)^{-1}(P-Q)(I+tP)^{-1} dt
\]
assuming that $B:=I+P>0$ and that $A:=I+Q>$. Rewrite this expression in terms of $A$ and $B$ 
and change the integration variable to $\tau=\frac{t-1}{t}$, to obtain
\[
\log(B)-\log(A)=\int_0^\infty (A+\tau I)^{-1}(B-A)(B+\tau I)^{-1}d\tau.
\]
If $B=A+\Delta$, then
\[(A+\tau I)^{-1}-(A+\Delta+\tau I)^{-1}= (A+\tau I)^{-1}\Delta (A+\Delta+\tau I)^{-1}
\]
which, for $A,\Delta$ Hermitian, $A>0$ and $A+\Delta>0$, leads to
\begin{eqnarray}\nonumber \log(A+\Delta)&=&\log(A)+\int_0^\infty (A+\tau I)^{-1}\Delta (A+\tau I)^{-1}d\tau\\
&&\hspace*{-65pt}+ \int_0^\infty (A+\tau I)^{-1}\Delta (A+\tau I)^{-1}\Delta (A+\tau I)^{-1}d\tau \label{loghessian}\\
&&\hspace*{-65pt}+ \;o(\|\Delta\|^2).\nonumber
\end{eqnarray}

By expanding in terms of eigenvectors of $A$ it can be verified that
\begin{equation}\label{eigenexpansion}
\int_0^\infty (A+\tau I)^{-1}\Delta (A+\tau I)^{-1}d\tau=M_A^{-1}(\Delta).
\end{equation}
Indeed, if $A$ is ${\rm diagonal}\{a_1,\ldots,a_n\}$ then
the $(i,j)$th entry of
\[
M_A\left( \int_0^\infty (A+\tau I)^{-1}\Delta (A+\tau I)^{-1} d\tau\right)
\]
is simply
\[
\int_0^1 a_i^{ (1-t)} a_j^t dt \int_0^\infty (a_i+\tau)^{-1}(a_j+\tau)^{-1}d\tau   (\Delta)_{ij} 
\]
where $(\Delta)_{ij} $ is the $(i,j)$th entry of $\Delta$ (in this same basis where $A$ is diagonal). Then
\[\int_0^1 a_i^ {(1-t)} a_j^t dt = \frac{a_i-a_j}{\log(a_i)-\log(a_j)}
\]
whereas, $\int_0^\infty (a_i+\tau)^{-1}(a_j+\tau)^{-1}d\tau$ is the inverse of the same expression.
This result is attributed to Lieb (see \cite[page 4]{ruskai2004}) and summarized below.

\begin{prop}\label{diffoflog}
{\sf The differential of $\log(A)$ is $M_A^{-1}$.}
\end{prop}

\end{document}